# REGENERATIVE TREE GROWTH: BINARY SELF-SIMILAR CONTINUUM RANDOM TREES AND POISSON–DIRICHLET COMPOSITIONS[1]

### By Jim Pitman and Matthias Winkel

#### *University of California, Berkeley and University of Oxford*

We use a natural ordered extension of the Chinese Restaurant Process to grow a two-parameter family of binary self-similar continuum fragmentation trees. We provide an explicit embedding of Ford's sequence of alpha model trees in the continuum tree which we identified in a previous article as a distributional scaling limit of Ford's trees. In general, the Markov branching trees induced by the two-parameter growth rule are not sampling consistent, so the existence of compact limiting trees cannot be deduced from previous work on the sampling consistent case. We develop here a new approach to establish such limits, based on regenerative interval partitions and the urn-model description of sampling from Dirichlet random distributions.

**1. Introduction.** We are interested in growth schemes for random rooted binary trees $T_n$ with $n$ leaves labeled by $[n] = \{1, \ldots, n\}$ of the following general form.

DEFINITION 1. Let $T_1$ be the tree with a single edge joining a root vertex and a leaf vertex labeled 1. Let $T_2$ be the Y-shaped tree consisting of a root and leaves labeled 1 and 2, each connected by an edge to a unique branch point.

To create $T_{n+1}$ from $T_n$, select an edge of $T_n$, say, $a_n \to c_n$, directed away from the root, replace it by three edges $a_n \to b_n$, $b_n \to c_n$ and $b_n \to n+1$ so that two new edges connect the two vertices $a_n$ and $c_n$ to a new branch point $b_n$ and a further edge connects $b_n$ to a new leaf labeled $n+1$.

Received March 2008.

[1]Supported in part by the NSF Awards 0405779 and 0806118, and by EPSRC Grant GR/T26368/01.

A *binary tree growth process* is a sequence $(T_n, n \geq 1)$ of random trees constructed in this way where at each step the edge $a_n \to c_n$ is selected randomly according to some *selection rule*, meaning a conditional distribution given $T_n$ for an edge of $T_n$. Given a selection rule, each tree $T_n$ has a distribution on the space $\mathbb{T}_{[n]}$ of rooted binary trees with $n$ leaves labeled $[n]$, and the selection rule specifies for all $n \geq 1$ conditional distributions of $T_{n+1}$ given $T_n$.

The *uniform rule*, where each of the $2n - 1$ edges of $T_n$ is selected with equal probability, gives a known binary tree growth process [25] related to the Brownian continuum random tree [1, 24]. Ford [10] introduced a one-parameter family of binary tree growth processes, where the selection rule for $0 < \alpha < 1$ is as follows:

(i) Given $T_n$ for $n \geq 2$, assign a weight $1 - \alpha$ to each of the $n$ edges adjacent to a leaf, and a weight $\alpha$ to each of the $n - 1$ other edges.
(ii) Select an edge of $T_n$ at random with probabilities proportional to the weights assigned by step (i).

For us, this selection rule will be the $(\alpha, 1 - \alpha)$-*rule*. Note that $\alpha = 1/2$ gives the uniform rule.

In [18] we showed that, also for $\alpha \neq 1/2$, the trees $T_n$ with leaf labels removed, denoted $T_n^\circ$, have a *continuum fragmentation tree* $\mathcal{T}^\alpha$ as their *distributional* scaling limit, when considered as $\mathbb{R}$-trees with unit edge lengths: $n^{-\alpha} T_n^\circ \to \mathcal{T}^\alpha$ in distribution for the Gromov–Hausdorff topology. However, in the main part of [18] and in all other fragmentation literature we are aware of, the labeling of leaves is *exchangeable*, while the labeling of leaves in order of appearance in the trees $T_n$ grown using the $(\alpha, 1 - \alpha)$-rule is not. Our results in [18] applied because of a weak sampling consistency of the $(\alpha, 1 - \alpha)$-trees; cf. [10]. The subtlety with these trees is that they are strongly sampling consistent in the sense defined in Definition 2 only if $\alpha = 1/2$; cf. [18].

DEFINITION 2. A binary tree growth process $(T_n, n \geq 1)$ is called *weakly sampling consistent* if the distributions of the delabeled trees $T_n^\circ$ and $\hat{T}_n^\circ$ coincide for all $n \geq 1$, where $\hat{T}_n^\circ$ is obtained from $T_{n+1}^\circ$ by removal of a leaf chosen uniformly at random. The process is called *strongly sampling consistent* if the distributions of $(T_n^\circ, T_{n+1}^\circ)$ and $(\hat{T}_n^\circ, T_{n+1}^\circ)$ coincide for all $n \geq 1$.

In this paper we take up the study of nonexchangeable labeling and the role of weak sampling consistency for a two-parameter extension of the $(\alpha, 1 - \alpha)$-rule; cf. Figure 1.

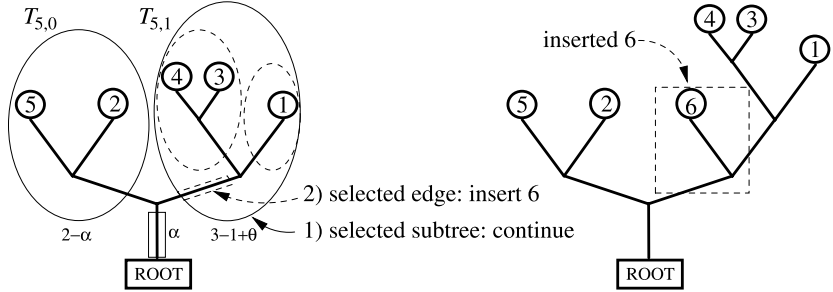

Fig. 1. *Recursive tree growth: in this scenario, the recursion consists of two steps. Weights for root edge and subtrees are displayed for the first step. The subtree $T_{n,1}$ is selected. Within tree $T_{n,1}$, the root edge is selected. Leaf 6 is inserted at the selected edge.*

DEFINITION 3. Let $0 \leq \alpha \leq 1$ and $\theta \geq 0$. We define the $(\alpha, \theta)$-*selection rule* as follows:

(i)$_{\mathrm{rec}}$ For $n \geq 2$, the tree $T_n$ branches at the branch point adjacent to the root into two subtrees $T_{n,0}$ and $T_{n,1}$. Given these are of sizes $m$ and $n - m$, say, where $T_{n,1}$ contains the smallest label in $T_n$, assign weight $\alpha$ to the edge connecting the root and the adjacent branch point, weights $m - \alpha$ and $n - m - 1 + \theta$, respectively, to the subtrees.

(ii)$_{\mathrm{rec}}$ Select the root edge or a subtree with probabilities proportional to these weights. If a subtree with two or more leaves was selected, recursively apply the weighting procedure (i)$_{\mathrm{rec}}$ to the selected subtree, until the root edge or a subtree with a single leaf was selected. If a subtree with a single leaf was selected, select the unique edge of this subtree.

A binary tree growth process $(T_n^{\alpha,\theta}, n \geq 1)$ grown via the $(\alpha, \theta)$-rules (i)$_{\mathrm{rec}}$, (ii)$_{\mathrm{rec}}$, for some $0 \leq \alpha \leq 1$ and $\theta \geq 0$, is called an $(\alpha, \theta)$-*tree growth process*.

For $\theta = 1 - \alpha$, each edge is chosen with the same probabilities as with Ford's rules (i) and (ii).

The boundary cases $\alpha = 0$ and $\alpha = 1$ are special and easy to describe (see Section 3.2). Growth is then linear or logarithmic in height, and scaling limits have a degenerate branching structure. We therefore focus on the parameter range $0 < \alpha < 1$ and study scaling limits and asymptotics of the associated trees $T_n = T_n^{\alpha,\theta}$.

We pointed out in [18] that Ford's $(\alpha, 1 - \alpha)$-tree growth process is associated with a *Chinese Restaurant Process (CRP)* as follows. The height $K_n$ of leaf 1 in $T_n$ increases whenever an edge on the path connecting 1 with the root, which we call the *spine*, is selected. Whenever a spinal edge is selected, the edge is replaced by two new spinal edges and a new subtree starts growing off the spine. If we call the subtrees off the spine *tables* and the leaves in subtrees *customers*, then the process of table sizes follows the $(\alpha, 1 - \alpha)$ seating plan of a CRP in the terminology of [24]. Similarly, we identify an $(\alpha, \theta)$ seating plan in the two-parameter model, meaning that the $(n+1)$st customer is seated at the $j$th table, with $n_j$ customers already seated, with probability proportional to $n_j - \alpha$ and at a new table with probability proportional to $\theta + k\alpha$, if $k$ tables are occupied. See Figure 2. Note that

$$\text{the } k\text{th customer in the restaurant is labeled } (k+1) \text{ as leaf in the tree,} \tag{1}$$

since leaf 1 is not in a subtree off the spine.

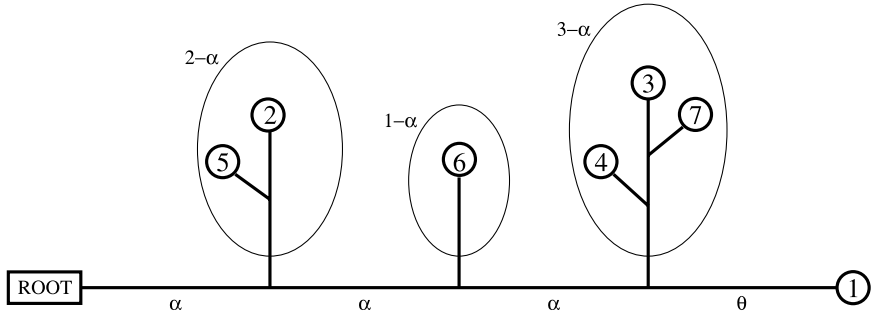

Fig. 2.  *The $(\alpha, \theta)$ tree growth procedure induces an ordered Chinese Restaurant Process.*

The theory of CRPs [24] immediately gives us a.s. a limit height $L_{\alpha, \theta} = \lim_{n \to \infty} K_n / n^\alpha$ of leaf 1, as well as limiting proportions $(P_1, P_2, \ldots)$ of leaves in each subtree *in birth order*, that is, in the order of least numbered leaves of subtrees, which can be represented as

$$(P_1, P_2, \ldots) = (W_1, \overline{W}_1 W_2, \overline{W}_1 \overline{W}_2 W_3, \ldots),$$

where the $W_i$ are independent, $W_i$ has a beta$(1 - \alpha, \theta + i\alpha)$ distribution on the unit interval, and $\overline{W}_i := 1 - W_i$. The distribution of the sequence of *ranked* limiting proportions is then Poisson–Dirichlet with parameters $(\alpha, \theta)$, for short PD$(\alpha, \theta)$.

However, this *spinal decomposition* of the tree also specifies the *spinal order*, that is, the order in which subtrees are encountered on the spine from the root to leaf 1 (from left to right in Figure 2). Note that due to the leaf labeling and the sequential growth of $T_n$, $n \geq 1$, subtrees are identifiable and keep their order throughout, which makes the spinal order consistent as $n$ varies. After the insertion of leaf $n + 1$, the sizes of subtrees in birth order and in spinal order form two compositions of $n$, $n \geq 1$. While the birth order is well known to be size-biased, we show that the compositions in spinal order form a regenerative composition structure in the sense of Gnedin and

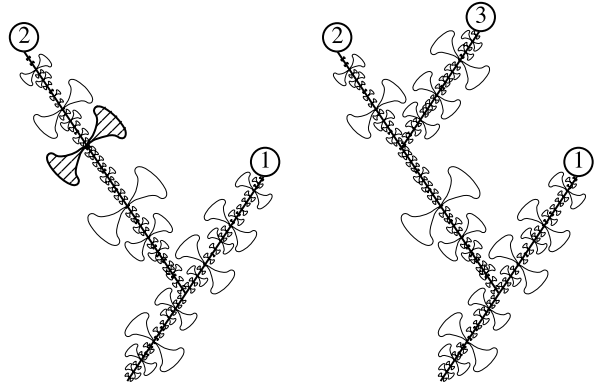

Fig. 3. *A tree equipped with strings of beads; crushing a bead into a new string of beads.*

Pitman [13], which is weakly sampling consistent for all $0 \leq \alpha \leq 1$ and $\theta \geq 0$, but not strongly so unless $\theta = \alpha$ [Proposition 6(i) and (ii)].

It follows from [13] in the strongly sampling consistent case $\theta = \alpha$ that the rescaled compositions converge almost surely to the associated regenerative interval partition and that the block containing leaf 2 is a size-biased pick from the composition of $n$, or from the interval partition in the limit $n \to \infty$. We obtain almost sure limiting results for the nonstrongly sampling consistent compositions (and discrete local times) in spinal order [Proposition 6(iii) and (iv)], and we solve the problem of finding leaf 2 in the nonstrongly consistent case, for the spinal composition of $n$ (Lemma 9) and for the limiting interval partition (Proposition 10). The limiting interval partition arranges the limiting proportions $(P_1, P_2, \ldots)$ in spinal order. We consider inverse local time $L^{-1}$ as a random distribution function on the interval $[0, L_{\alpha, \theta}]$. Then $([0, L_{\alpha, \theta}], dL^{-1})$ is an $(\alpha, \theta)$-string of beads in the following sense.

DEFINITION 4. An interval $(I, \mu)$ equipped with a discrete measure $\mu$ is called a *string of beads*. We refer to the weighted random interval $([0, L_{\alpha, \theta}], dL^{-1})$ associated with an $(\alpha, \theta)$-regenerative partition as $(\alpha, \theta)$-string of beads. We will also use this term for isometric copies of weighted intervals as in Figure 3.

As a by-product of these developments (Corollary 8), we obtain a sequential construction of the interval partition associated with the $(\alpha, \theta)$ regenerative composition structure described in [13], Section 8. This provides a much more combinatorial approach to the $(\alpha, \theta)$ regenerative interval partition than was given in [13], and solves the problem, left open in [13], of explicitly describing for general $(\alpha, \theta)$ how interval lengths governed by PD$(\alpha, \theta)$ should be ordered to form an $(\alpha, \theta)$ regenerative interval partition of $[0, 1]$ (Corollary 7).

We formulate and prove these results in Section 2. While they are key results for the study of the trees $T_n^{\alpha, \theta}$, they are also of independent interest in a framework of an *ordered CRP*. This notion will be made precise there and studied in some detail.

In Section 3 we formally introduce leaf-labeled rooted binary trees and the Markov branching property. We show that the delabeled trees from the $(\alpha, \theta)$-tree growth rules have the Markov branching property, and that the labeled trees have a regenerative property, which reflects the recursive nature of the growth rules (Proposition 11). We then study sampling consistency as defined in Definition 2:

PROPOSITION 1.  *Let* $(T_n^{\alpha, \theta}, n \geq 1)$ *be an* $(\alpha, \theta)$*-tree growth process for some* $0 < \alpha < 1$ *and* $\theta \geq 0$*, and* $T_n^{\alpha, \bar{\theta}, \circ}$*, $n \geq 1$, the associated delabeled trees.*

(a) $T_n^{\alpha, \theta}$ *has exchangeable leaf labels for all* $n \geq 1$ *if and only if* $\alpha = \theta = 1/2$.
(b) $(T_n^{\alpha, \theta, \circ}, n \geq 1)$ *is strongly sampling consistent if and only if* $\alpha = \theta = 1/2$.
(c) $(T_n^{\alpha, \theta, \circ}, n \geq 1)$ *is weakly sampling consistent if and only if* $\theta = 1 - \alpha$ *or* $\theta = 2 - \alpha$.

We actually show that the distributions of delabeled trees coincide for $\theta = 1 - \alpha$ and $\theta = 2 - \alpha$, and do so only in these weakly sampling consistent cases (Lemma 12).

The main contribution of this paper is to establish limiting continuum random trees (CRTs) even without weak sampling consistency. For a tree $T_n$ labeled by $[n] = \{1, \ldots, n\}$, we denote by $S(T_n; [k])$ the smallest subtree of $T_n$ that contains the root and the leaves labeled $1, \ldots, k$. It will be convenient to use Aldous's formalism of reduced trees with edge lengths: denote by $R(T_n; [k])$ the tree $T_k$ with edges marked as follows; because of the growth procedure each vertex of $T_k$ is also a vertex of $T_n$, and we mark each edge of $T_k$ by the graph distance in $T_n$ of the two vertices that the edge connects. First, we study the asymptotics of these reduced trees.

PROPOSITION 2.  *Let* $(T_n^{\alpha, \theta}, n \geq 1)$ *be an* $(\alpha, \theta)$*-tree growth process. If* $0 < \alpha < 1$ *and* $\theta \geq 0$*, then*

$$n^{-\alpha} R(T_n^{\alpha, \theta}, [k]) \to \mathcal{R}_k^{\alpha, \theta} \qquad \text{almost surely as } n \to \infty,$$

*in the sense that the* $2k - 1$ *edge lengths of* $R(T_n^{\alpha, \theta}, [k])$ *scaled by* $n^\alpha$ *converge almost surely as* $n \to \infty$ *to limiting edge lengths of a tree* $\mathcal{R}_k^{\alpha, \theta}$*, for all* $k \geq 1$.

We proved this in [18], Proposition 18, for Ford's $(\alpha, 1 - \alpha)$-tree growth process. As in [18], we will also provide an explicit description of the distribution of $(\mathcal{R}_k^{\alpha, \theta}, k \geq 1)$. We will, in fact, prove a stronger statement for trees $\mathcal{R}_k^{\alpha,\theta}$ where each edge has the structure of a string of beads that records limiting proportions of leaves of subtrees as atoms on the branches (Proposition 14 and Corollary 15). We deduce growth rules for the passage from $k$ to $k+1$ leaves for the limiting trees equipped with strings of beads (Corollary 16). These are remarkably simple and consist of picking a bead (using Proposition 10) and crushing the bead of size $s_k$, say, into $m_{k+1}$, where $m_{k+1}/s_k \sim \mathrm{PD}(\alpha,\theta)$, arranging these as a new string of beads (using Corollary 7), attaching them to the location of the bead, which now splits an edge and the remainder of its string of beads into two, as illustrated in Figure 3.

In the $(\alpha, 1-\alpha)$ case, growth by crushing beads is closely connected to growth rules for random recursive trees studied by Dong, Goldschmidt and Martin [6]. Specifically, we can associate with $\mathcal{R}_k$ a tree $V_k$ with $k$ *vertices* labeled by $[k]$ and infinitely many unlabeled vertices, all marked by weights; let $V_1$ consist of a root labeled 1 and infinitely many unlabeled children marked by the sequence $m_1$ of masses of the string of beads on $\mathcal{R}_1$; to construct $V_{k+1}$ from $V_k$, identify the unlabeled leaf in $V_k$ marked by the size of the chosen bead $s_k$, label it by $k+1$ and add infinitely many children of vertex $k+1$, marked by the sizes $m_{k+1}$ of the crushed bead. The limit $V_\infty$ is a recursive tree where all vertices have infinitely many children. We show in this paper that the richer structure of $(\mathcal{R}_k, \mu_k)$, that includes edges on which the atoms of $\mu_k$ are distributed, has a binary CRT as its limit. In fact, $V_\infty$ can be constructed for general $(\alpha, \theta)$, but the purpose of [6] was to establish a coagulation-fragmentation duality that only works for $(\alpha, 1-\alpha)$. See also Blei, Griffiths and Jordan [5] for another application of nested Chinese restaurant processes to define distributions on infinitely-deep, infinitely-branching trees.

Section 4 will establish CRT limits for the general $(\alpha, \theta)$-tree growth process.

THEOREM 3. *In the setting of Proposition 2, there exists a CRT $\mathcal{T}^{\alpha,\theta}$ on the same probability space such that we have for the delabeled trees $\mathcal{R}_k^{\alpha,\theta,\circ}$, $k \geq 1$, associated with $\mathcal{R}_k^{\alpha,\theta}$, $k \geq 1$, that*

$$\mathcal{R}_k^{\alpha,\theta,\circ} \to \mathcal{T}^{\alpha,\theta} \qquad \textit{almost surely as } k \to \infty, \textit{ in the Gromov–Hausdorff topology.}$$

In fact, CRTs such as $\mathcal{T}^{\alpha,\theta}$ are equipped with a mass measure $\mu$. We can construct $\mu$ as the limit of the strings of beads that we constructed on $\mathcal{R}_k^{\alpha,\theta,\circ}$ [see Corollary 23], using Evans' and Winter's [9] weighted Gromov–Hausdorff convergence that we recall in Section 4.1.

It would be nice to replace the two-step limiting procedure of Proposition 2 and Theorem 3 for trees reduced to $k$ leaves, letting first $n \to \infty$ and then $k \to \infty$, by a single statement:

CONJECTURE 1. *In the setting of Proposition 2, we have*

$$n^{-\alpha} \mathcal{T}_n^{\alpha,\theta,\circ} \to \mathcal{T}^{\alpha,\theta} \qquad \text{almost surely, as } n \to \infty,$$

*for the Gromov–Hausdorff topology.*

In [18] we used exchangeability to obtain fine tightness estimates and establish convergence in probability for a wide class of exchangeable strongly sampling consistent Markov branching trees. From this result we deduce a convergence in distribution in the weakly sampling consistent cases $\theta = 1 - \alpha$ and $\theta = 2 - \alpha$, but without sampling consistency, this argument breaks down.

Our method to prove Theorem 3 uses an embedding of $(T_n^{\alpha,\theta}, n \geq 1)$ and $(\mathcal{R}_k^{\alpha,\theta}, k \geq 1)$ in a given fragmentation CRT. For a rooted $\mathbb{R}$-tree $(\mathcal{T}, \rho)$ and leaves $\Sigma_1, \ldots, \Sigma_k$ of $\mathcal{T}$, denote by $R(\mathcal{T}; \Sigma_1, \ldots, \Sigma_k)$ the smallest subtree of $\mathcal{T}$ that contains $\rho$ and $\Sigma_1, \ldots, \Sigma_k$. The family of binary fragmentation CRTs $\mathcal{T}$ is parameterized by a self-similarity parameter $\alpha > 0$ and a dislocation measure $\nu(du)$, a sigma-finite measure on $[1/2, 1)$ with $\int_{[1/2,1)}(1-u)\nu(du) < \infty$; see Section 4.1.

THEOREM 4. *Let $(\mathcal{T}^{\alpha,\theta}, \rho, \mu)$ be a binary fragmentation CRT with root $\rho$ and mass distribution $\mu$, associated with dislocation measure $\nu_{\alpha,\theta}(du) = f_{\alpha,\theta}^\circ(u)\,du$, $1/2 < u < 1$, where*

$$\Gamma(1-\alpha)f_{\alpha,\theta}^\circ(u) = \alpha(u^\theta(1-u)^{-\alpha-1} + u^{-\alpha-1}(1-u)^\theta)$$
$$+ \theta(u^{\theta-1}(1-u)^{-\alpha} + u^{-\alpha}(1-u)^{\theta-1})$$

*for some $0 < \alpha < 1$ and $\theta \geq 0$. Then there exists, on a suitably extended probability space, a sequence $(\Sigma_n, n \geq 1)$ of random leaves of $\mathcal{T}^{\alpha,\theta}$, such that $(R(\mathcal{T}^{\alpha,\theta}; \Sigma_1, \ldots, \Sigma_k), k \geq 1)$ has the same distribution as $(\mathcal{R}_k^{\alpha,\theta}, k \geq 1)$.*

With this embedding, the projection of the mass distribution $\mu$ of $\mathcal{T}^{\alpha,\theta}$ onto $R(\mathcal{T}^{\alpha,\theta}; \Sigma_1, \ldots, \Sigma_k)$ yields strings of beads with distributions as we constructed them on $\mathcal{R}_k^{\alpha,\theta}$. See Proposition 20.

## 2. An ordered Chinese Restaurant Process and regenerative composition structures.

2.1. *Regenerative compositions.* We recall in this subsection some background on regenerative composition structures from [13]. A *composition of $n$* is a sequence $(n_1, \ldots, n_k)$ of positive integers with sum $n$. A sequence of random compositions $\mathcal{C}_n$ of $n$ is called *regenerative* if, conditionally given that the first part of $\mathcal{C}_n$ equals $n_1$, the remaining parts of $\mathcal{C}_n$ define a composition of $n - n_1$ with the same distribution as $\mathcal{C}_{n-n_1}$. Given any decrement matrix $(q(n, m), 1 \leq m \leq n)$, there is an associated sequence $\mathcal{C}_n$ of regenerative random compositions of $n$ defined by specifying that $q(n, \cdot)$ is the distribution of the first part of $\mathcal{C}_n$. Thus, for each composition $(n_1, \ldots, n_k)$ of $n$,

$$
\begin{aligned}
(2) \quad & \mathbb{P}(\mathcal{C}_n = (n_1, \ldots, n_k)) \\
& = q(n, n_1) q(n - n_1, n_2) \cdots q(n_{k-1} + n_k, n_{k-1}) q(n_k, n_k).
\end{aligned}
$$

We regard a composition of $n$ as a distribution of identical balls in an ordered sequence of boxes. For a sequence of compositions $(\mathcal{C}_n, n \geq 1)$, let $\hat{\mathcal{C}}_n$ denote the composition of $n$ obtained by removal of a ball chosen uniformly at random from $\mathcal{C}_{n+1}$, and discarding the empty box if the chosen ball is the only one in its box. We call $(\mathcal{C}_n, n \geq 1)$ *weakly sampling consistent* if $\mathcal{C}_n \overset{d}{=} \hat{\mathcal{C}}_n$ for every $n$, and *strongly sampling consistent* if $(\mathcal{C}_n, \mathcal{C}_{n+1}) \overset{d}{=} (\hat{\mathcal{C}}_n, \mathcal{C}_{n+1})$ for every $n$. A detailed theory of the asymptotic behavior of weakly sampling consistent sequences of regenerative compositions of $n$ (known as *composition structures*) is provided in [13].

Now write

$$
\mathcal{C}_n = (N_{n,1}, N_{n,2}, \ldots, N_{n,K_n}) \quad \text{and let} \quad S_{n,k} = \sum_{j=1}^{k} N_{n,j},
$$

where $N_{n,j} = 0$ for $j > K_n$. According to Gnedin and Pitman [13], if $(\mathcal{C}_n, n \geq 1)$ is weakly sampling consistent, there is the following convergence in distribution of random sets with respect to the Hausdorff metric on closed subsets of $[0, 1]$:

$$
(3) \quad \{S_{n,k}/n, k \geq 0\} \xrightarrow[n \to \infty]{\mathrm{d}} \mathcal{Z} := \{1 - \exp(-\xi_t), t \geq 0\}^{\mathrm{cl}},
$$

where the left-hand side is the random discrete set of values $S_{n,k}$ rescaled onto $[0, 1]$, and the right-hand side is the closure of the range of 1 minus the exponential of a subordinator $(\xi_t, t \geq 0)$. If $(\mathcal{C}_n, n \geq 1)$ is strongly sampling consistent, then the convergence (3) holds also with convergence in distribution replaced by almost sure convergence. The collection of open interval components of $[0, 1] \setminus \mathcal{Z}$ is then called the *regenerative interval partition* associated with $(\mathcal{C}_n, n \geq 1)$. In particular, a strongly sampling consistent composition structure can be derived from $\mathcal{Z}$ by uniform sampling in $[0, 1]$ using $\mathcal{Z}$ to separate parts.

The distribution of a subordinator $(\xi_t, t \geq 0)$ is encoded in its Laplace exponent $\Phi$ as

$$
\mathbb{E}(e^{-s\xi_t}) = e^{-t\Phi(s)} \quad \text{where } \Phi(s) = a + cs + \int_{(0,\infty)} (1 - e^{-sx}) \Lambda(dx),
$$

for all $s \geq 0$, $t \geq 0$, and characteristics $(a, c, \Lambda)$, where $a \geq 0$, $c \geq 0$ and $\Lambda$ is a measure on $(0, \infty)$ with $\int_{(0,\infty)} (1 \wedge x) \Lambda(dx) < \infty$.

2.2. *An ordered Chinese Restaurant Process.* We will now use an ordered version of the CRP to construct an exchangeable random partition $\Pi_{\alpha,\theta}$ of $\mathbb{N}$ governed by the CRP as described in [24], jointly with a random total ordering of the blocks (tables) of $\Pi_{\alpha,\theta}$. With a suitable encoding that we make precise, this random total ordering is independent of $\Pi_{\alpha,\theta}$.

First recall the $(\alpha,\theta)$ CRP for fixed $0 \le \alpha \le 1$ and $\theta \ge 0$. Customers labeled by $\mathbb{N} := \{1, 2, \ldots\}$ seat themselves at tables labeled by $\mathbb{N}$ in the following way: Customer 1 sits at table 1. Given that $n$ customers have been seated at $k$ different tables, with $n_i$ customers at table $i$ for $i \in [k]$, customer $n+1$

- sits at the $i$th occupied table with probability $(n_i - \alpha)/(n + \theta)$, for $i \in [k]$;
- sits alone at table $k + 1$ with probability $(k\alpha + \theta)/(n + \theta)$.

The state of the system after $n$ customers have been seated is a random partition $\Pi_n$ of $[n]$. By construction, these partitions are exchangeable, and consistent as $n$ varies so they induce a random partition $\Pi_\infty$ of $\mathbb{N}$ whose restriction to $[n]$ is $\Pi_n$.

When $\alpha = 1$, $\Pi_\infty$ consists of all singleton blocks since no customer ever sits at an occupied table. So we assume henceforth that $0 \le \alpha < 1$. Basic facts are that the block of $\Pi_\infty$ associated with table $j$ has an almost sure limiting frequency $P_j$, and that the $P_j$ may be represented as

$$(4) \qquad (P_1, P_2, \ldots) = (W_1, \overline{W}_1 W_2, \overline{W}_1 \overline{W}_2 W_3, \ldots),$$

where the $W_i$ are independent, $W_i \sim \text{beta}(1 - \alpha, \theta + i\alpha)$ and $\overline{W}_i := 1 - W_i$. Note that the proportions $(P_1, P_2, \ldots)$ are in a size-biased random order, corresponding to the fact that the table numbers label the blocks of $\Pi_\infty$ in order of their least elements.

Another basic fact, read from [24], is that the number $K_n$ of occupied tables after $n$ customers (number of blocks of $\Pi_n$) has the limiting behavior

$$(5) \qquad K_n/n^\alpha \xrightarrow{\text{a.s.}} L_{\alpha,\theta} = \Gamma(1 - \alpha) \lim_{j \to \infty} j(P_j^\downarrow)^\alpha \qquad \text{for } 0 < \alpha < 1,$$

where $(P_j^\downarrow, j \ge 1)$ is the *ranked* sequence of proportions $(P_j, j \ge 1)$, and $L_{\alpha,\theta}$ is a random variable with the tilted Mittag–Leffler distribution with moments

$$(6) \qquad \mathbb{E}(L_{\alpha,\theta}^n) = \frac{\Gamma(\theta + 1)}{\Gamma(\theta/\alpha + 1)} \frac{\Gamma(\theta/\alpha + n + 1)}{\Gamma(\theta + n\alpha + 1)} \qquad (n \ge 0).$$

This $L_{\alpha,\theta}$ is the local time variable associated with a regenerative $\text{PD}(\alpha, \theta)$ interval partition of $[0, 1]$, also called its $\alpha$-diversity. For $\alpha = 0$, we have $K_n/\log(n) \to \theta$ almost surely.

We now put a random total order $<$ on the tables as follows. Independently of the process of seating of customers at tables, let the tables be ordered from left to right according to the following scheme. Put the second table to the right of the first with probability $\theta/(\alpha + \theta)$ and to the left with probability $\alpha/(\alpha + \theta)$. This creates three possible locations for the third table: put it

- to the left of the first two tables with probability $\alpha/(2\alpha + \theta)$;
- between the first two tables with probability $\alpha/(2\alpha + \theta)$;
- to the right of the first two tables with probability $\theta/(2\alpha + \theta)$.

And so on: given any one of $k!$ possible orderings of $k$ tables from left to right, there are $k + 1$ possible places for the $(k + 1)$st table to be squeezed in. The place to the right of all $k$ tables is assigned probability $\theta/(k\alpha + \theta)$; each of the other $k$ places is assigned probability $\alpha/(k\alpha + \theta)$.

Let $\sigma_k(i)$ denote the location of the $i$th table relative to the first $k$ tables, counting from 1 for the left-most table to $k$ for the right-most. So $\sigma_k$ is a random permutation of $[k]$. The sequence of permutations $(\sigma_k, k \geq 1)$ is consistent in the sense that if $\sigma_k(i) < \sigma_k(j)$ for some $k \geq i \vee j := \max\{i, j\}$, then the same is true for all $k \geq i \vee j$. Thus, the sequence $(\sigma_k, k \geq 1)$ specifies a random total order on $\mathbb{N}$, call it the *table order*. Given $\sigma_1, \ldots, \sigma_k$,

- $\sigma_{k+1}(k + 1) = k + 1$ with probability $\theta/(k\alpha + \theta)$;
- $\sigma_{k+1}(k + 1) = i$ with probability $\alpha/(k\alpha + \theta)$ for each $i \in [k]$

and

$$(7) \qquad \sigma_{k+1}(j) = \sigma_k(j) + 1(\sigma_{k+1}(k + 1) \leq \sigma_k(j)) \qquad \text{for } j \in [k].$$

Thus, by construction, $(\sigma_k, k \geq 1)$ is independent of the (unordered) random partition $\Pi_\infty$ of $\mathbb{N}$, with

$$\mathbb{P}(\sigma_k = \pi) = \frac{(\theta/\alpha)^{R(\pi)}}{[\theta/\alpha]_k}$$

for each permutation $\pi$ of $[k]$, where

$$[x]_k = x(x + 1)\cdots(x + k - 1) = \Gamma(x + k)/\Gamma(x)$$

and

$$R(\pi) := \sum_{i=1}^{k} 1(\pi_j > \pi_i \text{ for all } 1 \leq i < j)$$

is the number of record values in the permutation $\pi$. Note that for $k \geq 2$ the distribution of $\sigma_k$ is uniform iff $\alpha = \theta$. The formulas apply as suitable limit expressions: if $\alpha = 0$ and $\theta > 0$, tables are ordered in order of appearance and $\sigma_k$ is the identity permutation (there is only one table for $\theta = 0$); if $0 < \alpha < 1$ and $\theta = 0$, the first table remains right-most, and the $\sigma_k$ is uniform on permutations with $\pi(1) = k$. See [21, 23] for related work.

In the sequel, we will repeatedly use generalized urn scheme arguments, so let us briefly review the main points here. See [22] and [24], Section 2.2, for references. Recall that the distribution of a random vector $\Delta = (\Delta_1, \ldots, \Delta_m)$ with $\Delta_1 + \cdots + \Delta_m = 1$ and density

$$g_{\Delta_1, \ldots, \Delta_{m-1}}(x_1, \ldots, x_{m-1})$$
$$= \frac{\Gamma(\gamma_1 + \cdots + \gamma_m)}{\Gamma(\gamma_1) \cdots \Gamma(\gamma_m)} x_1^{\gamma_1 - 1} \cdots x_{m-1}^{\gamma_{m-1} - 1} (1 - x_1 - \cdots - x_{m-1})^{\gamma_m - 1}$$

on $\{(x_1, \ldots, x_{m-1}) : x_1, \ldots, x_{m-1} \geq 0, x_1 + \cdots + x_{m-1} \leq 1\}$ is called the *Dirichlet distribution* with parameters $\gamma_1, \ldots, \gamma_m > 0$.

LEMMA 5. (i) *Consider a weight vector* $\gamma = (\gamma_1, \ldots, \gamma_m)$ *and a process* $(H^{(n)}, n \geq 0)$ *with* $H^{(0)} = 0$, *where* $H^{(n)} = (H_1^{(n)}, \ldots, H_m^{(n)})$ *evolves according to the updating rule to increase by 1 a component chosen with probabilities proportional to current weights* $\gamma + H^{(n)}$:

$$\mathbb{P}(H^{(n+1)} = H^{(n)} + e_i \mid H^{(1)}, \ldots, H^{(n)})$$
$$= \frac{\gamma_i + H_i^{(n)}}{\gamma_1 + \cdots + \gamma_m + n} \qquad a.s., \; i = 1, \ldots, m,$$

*where* $e_i$ *is the* $i$th *unit vector. Then* $H^{(n)}/n \xrightarrow[n \to \infty]{\text{a.s.}} \Delta \sim \mathrm{Dirichlet}(\gamma_1, \ldots, \gamma_m)$ *and*

$$\mathbb{P}(H^{(n+1)} = H^{(n)} + e_i \mid H^{(1)}, \ldots, H^{(n)}, \Delta) = \Delta_i \qquad a.s., \; i = 1, \ldots, m,$$

*which means that the components of increase are conditionally independent and identically distributed according to the limiting weight proportions* $\Delta$.

(ii) *A vector* $\Delta \sim \mathrm{Dirichlet}(\gamma_1, \ldots, \gamma_m)$ *can be represented as*

$$(\Delta_1, \ldots, \Delta_m) = (W_1, \overline{W_1} W_2, \overline{W_1}\,\overline{W_2} W_3, \ldots,$$
$$\overline{W_1} \cdots \overline{W_{m-2}} W_{m-1}, \overline{W_1} \cdots \overline{W_{m-2}}\,\overline{W_{m-1}}),$$

*where the* $W_i$ *are independent,* $W_i \sim \mathrm{beta}(\gamma_i, \gamma_{i+1} + \cdots + \gamma_m)$ *and* $\overline{W_i} := 1 - W_i$.

If $\gamma \in \mathbb{N}^m$, the process $H$ arises when drawing from an urn with initially $\gamma_i$ balls of color $i$, always adding a ball of the color drawn.

2.3. *The composition of table sizes in the ordered Chinese Restaurant Process.* Let $\widetilde{\Pi}_n$ denote the random ordered partition of $[n]$ induced by ordering the blocks of $\Pi_n$ according to $\sigma_{K_n}$, where $K_n$ is the number of blocks of $\Pi_n$. Let $\mathcal{C}_n$ denote the random composition of $n$ defined by the sizes of blocks of $\widetilde{\Pi}_n$. If $\mathcal{C}_n^*$ is the sequence of sizes of blocks of $\Pi_n$, in order of least elements (or table label), and $K_n = k$, the $j$th term of $\mathcal{C}_n^*$ is the $\sigma_k(j)$th term of $\mathcal{C}_n$.

PROPOSITION 6. (i) *For each $(\alpha, \theta)$ with $0 < \alpha < 1$ and $\theta \geq 0$ the sequence of compositions $(\mathcal{C}_n, n \geq 1)$ defined as above is regenerative, with decrement matrix*

$$(8) \quad q_{\alpha, \theta}(n, m) = \binom{n}{m} \frac{n\alpha - m\alpha + m\theta}{n} \frac{[1 - \alpha]_{m-1}}{[n - m + \theta]_m} \qquad (1 \leq m \leq n).$$

(ii) *This sequence of compositions $(\mathcal{C}_n, n \geq 1)$ is weakly sampling consistent, but strongly sampling consistent only if $\alpha = \theta$.*

(iii) *Let $S_{n,j}$ be the number of the first $n$ customers seated in the $j$ leftmost tables. Then there is the following almost sure convergence of random sets with respect to the Hausdorff metric on closed subsets of $[0, 1]$:*

$$(9) \qquad \{S_{n,j}/n, j \geq 0\} \xrightarrow[n \to \infty]{a.s.} \mathcal{Z}_{\alpha, \theta} := \{1 - \exp(-\xi_t), t \geq 0\}^{\mathrm{cl}},$$

*where the left-hand side is the random discrete set of values $S_{n,j}$ rescaled onto $[0, 1]$, and the right-hand side $\mathcal{Z}_{\alpha, \theta}$ is by definition the closure of the range of 1 minus the exponential of the subordinator $(\xi_t, t \geq 0)$ with Laplace exponent*

$$(10) \quad \begin{aligned} \Phi_{\alpha, \theta}(s) &= \frac{s\Gamma(s + \theta)\Gamma(1 - \alpha)}{\Gamma(s + \theta + 1 - \alpha)} \qquad for\ \theta > 0 \quad and \\ \Phi_{\alpha, 0}(s) &= \frac{\Gamma(s + 1)\Gamma(1 - \alpha)}{\Gamma(s + 1 - \alpha)}. \end{aligned}$$

(iv) *Also, if $L_n(u)$ denotes the number of $j \in \{1, \ldots, K_n\}$ with $S_{n,j}/n \leq u$, then*

$$(11) \qquad \lim_{n \to \infty} \sup_{u \in [0,1]} |n^{-\alpha} L_n(u) - L(u)| = 0 \qquad a.s.,$$

*where $L := (L(u), u \in [0, 1])$ is a continuous local time process for $\mathcal{Z}_{\alpha, \theta}$, meaning that the random set of points of increase of $L$ is $\mathcal{Z}_{\alpha, \theta}$ almost surely.*

NOTE. Various characterizations of $L$ can be given in terms of $\mathcal{Z}_{\alpha, \theta}$ and $\xi$. See below.

PROOF OF PROPOSITION 6. (i) That $(\mathcal{C}_n)$ is regenerative is proved by induction on $n$. The case $n = 1$ is trivial, and if $(\mathcal{C}_m, 1 \leq m \leq n)$ is regenerative, then, by the seating rule, three scenarios can occur. Given customer $n + 1$ sits alone at a new first table, the remaining composition $\mathcal{C}_n$ is trivially distributed as $\mathcal{C}_n$. Given customer $n + 1$ sits down at the existing first table of size $n_1$, the induction hypothesis implies that the remaining composition is distributed as $\mathcal{C}_{n-n_1}$, as required. Given customer $n + 1$ sits neither at a new first nor at the existing first table of size $n_1$, the seating rules are such that he chooses his seat in the remaining composition as if he were customer $n - n_1 + 1$ for composition $\mathcal{C}_{n-n_1}$, and the induction hypothesis allows to conclude that the resulting composition of $n - n_1 + 1$ is distributed as $\mathcal{C}_{n-n_1+1}$, as required.

Denote by $q(n, m)$ the probability that the first block in $\mathcal{C}_n$ is of size $m$. Then, the seating rules imply that

$$q(n+1, m) = q(n, m-1)\frac{m-1-\alpha}{n+\theta} + q(n, m)\frac{n+\theta-m}{n+\theta}$$

$$(12) \qquad\qquad + \frac{\alpha}{n+\theta}1_{\{m=1\}}, \qquad 1 \le m \le n+1,$$

where $q(n, m) = 0$ for $m > n$ or $m \le 0$. It is enough to check that the matrix given in (8) solves (12) for $m \ge 2$, that is to show

$$\binom{n+1}{m}\frac{n\alpha + \alpha - m\alpha + m\theta}{n+1}\frac{[1-\alpha]_{m-1}}{[n+1-m+\theta]_m}$$

$$= \binom{n}{m-1}\frac{n\alpha - m\alpha + \alpha + m\theta - \theta}{n}\frac{[1-\alpha]_{m-2}}{[n-m+1+\theta]_{m-1}}\frac{m-1-\alpha}{n+\theta}$$

$$+ \binom{n}{m}\frac{n\alpha - m\alpha + m\theta}{n}\frac{[1-\alpha]_{m-1}}{[n-m+\theta]_m}\frac{n+\theta-m}{n+\theta}.$$

Obvious cancellations reduce this to

$$n(n\alpha + \alpha - m\alpha + m\theta) = m(n\alpha - m\alpha + \alpha + m\theta - \theta)$$

$$+ (n+1-m)(n\alpha - m\alpha + m\theta),$$

which is easily verified. The decrement matrix (8) was derived in [13], Section 8, as that associated with the unique regenerative composition structure whose interval partition of $[0, 1]$ has ranked lengths distributed according to the Poisson–Dirichlet distribution with parameters $(\alpha, \theta)$. Thus, formula (8) gives the decrement matrix of a weakly sampling consistent family of regenerative compositions.

(ii) Weak sampling consistency was a by-product of the proof of (i). Let us show that $(\mathcal{C}_n, n \ge 1)$ is strongly sampling consistent if and only if $\alpha = \theta$. It is known that the compositions induced by independent uniform variables separated by the zero-set of a $(2 - 2\alpha)$-dimensional Bessel bridge have the dynamics of the Chinese Restaurant Process with seating plan $(\alpha, \alpha)$ and a uniform block order. Also, this construction using a Bessel bridge generates a strongly sampling consistent composition structure. On the other hand, the ordered version of the Chinese Restaurant Process also induces a uniform block order for $\alpha = \theta$. Conversely, calculate the following probabilities:

$$\mathbb{P}(\mathcal{C}_2 = (1, 1)) = \frac{\alpha + \theta}{1 + \theta}, \qquad \mathbb{P}(\mathcal{C}_2 = (2)) = \frac{1 - \alpha}{1 + \theta},$$

$$\mathbb{P}(\mathcal{C}_3 = (2, 1)) = \frac{(\alpha + 2\theta)(1 - \alpha)}{(1 + \theta)(2 + \theta)},$$

and note that strong sampling consistency requires

$$\frac{(1-\alpha)\theta}{(1+\theta)(2+\theta)} = \mathbb{P}(\mathcal{C}_2 = (2), \mathcal{C}_3 = (2,1)) = \frac{(\alpha+2\theta)(1-\alpha)}{(1+\theta)(2+\theta)}\frac{1}{3}$$

$$\iff \quad \alpha = \theta.$$

(iii) Now (3) yields convergence in distribution in (9), and (10) was derived in [13], formula (41). To get the almost sure convergence in (9), observe that for each $i \geq 1$, the proportion $P_i^{(n)}$ of customers at the $i$th table in order of appearance corresponds to the size of a gap in $\{S_{n,j}/n, j \geq 1\}$ and converges to $P_i$ almost surely as $n \to \infty$. As for the gap $(G_i^{(n)}, D_i^{(n)})$ itself, where $D_i^{(n)} = G_i^{(n)} + P_i^{(n)}$, a simple argument allows to also deduce almost sure convergence as $n \to \infty$,

$$G_i^{(n)} = \frac{S_{n,\sigma_{K_n}(i)-1}}{n} = \sum_{j=1}^{\infty} P_j^{(n)} 1_{\{\sigma_{j\vee i}(j) < \sigma_{j\vee i}(i)\}}$$

$$\to \sum_{j=1}^{\infty} P_j 1_{\{\sigma_{j\vee i}(j) < \sigma_{j\vee i}(i)\}} =: G_i,$$

and, hence, $D_i^{(n)} \to G_i + P_i =: D_i$, using the consistent construction of the sequence $(\sigma_k, k \geq 1)$ and the almost sure convergence of frequencies of all classes of $\Pi_\infty$.

In particular, on a set of probability one, the following holds. For each $\varepsilon > 0$ the locations of all gaps of length $P_i > \varepsilon$ converge, and a simple argument shows that we can find $n_0 \geq 1$ such that, for all $n \geq n_0$,

$$B(\{S_{n,j}/n, j \geq 1\}, \varepsilon) \supset \{G_i, D_i, i \geq 1\} \quad \text{and}$$

$$B(\{G_i, D_i, i \geq 1\}, \varepsilon) \supset \{S_{n,j}/n, j \geq 1\},$$

where $B(\mathcal{S}, \varepsilon) = \{x \in [0,1] : |x - y| \leq \varepsilon \text{ for some } y \in \mathcal{S}\}$ for any Borel set $\mathcal{S} \subset [0,1]$. We deduce the almost sure Hausdorff convergence of (9). Cf. the arXiv version [11] of [12] for a similar argument.

(iv) As for convergence of local time processes, the convergence (5) of $L_n(1)/n^\alpha = K_n/n^\alpha$ to $L(1)$ equal to the $\alpha$-diversity of the limiting $\mathrm{PD}(\alpha, \theta)$ is established in [24]. Look next at a time $u$ in the random interval $(G_1, D_1)$ associated with the first table. The dynamics of the table ordering imply that the numbers of tables to the left of the first table develop according to the urn scheme associated with sampling from a beta$(1, \theta/\alpha)$ variable $\beta_{1,\theta/\alpha}$ which is independent of $L(1)$. It follows that for $u$ in $(G_1, D_1)$ there is almost sure convergence of $L_n(u)/n^\alpha$ to $\beta_{1,\theta/\alpha}L(1)$. Similarly, if we look at the first $k$ tables, and count how numbers of following tables fall in the $k+1$ gaps they create, we see the dynamics associated with sampling from a Dirichlet distribution with its first $k$ parameters equal to 1 and the last equal to $\theta/\alpha$; cf. Lemma 5. As $k \to \infty$, the associated cumulative Dirichlet fractions are almost surely dense in $[0, 1]$. It follows that we get a.s. convergence in (11) for all $u$ in the random set of times $\bigcup_{j \geq 1}(G_j, D_j)$, and that the countable random set of a.s. distinct limit values from these intervals is a.s. dense in $[0, L(1)]$. The conclusion then follows by a standard argument; cf. [15]. $\quad\square$

It is worth recording some consequences of this argument.

COROLLARY 7. *The collection of intervals*

$$\bigcup_{j \geq 1}(G_j, D_j)$$

*for $(G_j, D_j, j \geq 1)$ created from the size-biased frequencies $(P_j, j \geq 1)$ and the independent sequence of random permutations $(\sigma_k, k \geq 1)$ specified in (7) provides an explicit construction of a regenerative $(\alpha, \theta)$ interval partition of $[0, 1]$.*

COROLLARY 8. *Construct a random interval partition of $[0, 1]$ as follows. Let $(G_1, D_1)$ be such that the joint law of $(G_1, D_1 - G_1, 1 - D_1)$ is Dirichlet$(\alpha, 1 - \alpha, \theta)$ for some $0 < \alpha < 1$ and $\theta \geq 0$. Given $(G_1, D_1)$, let this be one interval component, let the interval components within $[0, G_1]$ be obtained by linear scaling of a regenerative $(\alpha, \alpha)$ partition, and let the interval components within $[D_1, 1]$ be obtained by linear scaling of a regenerative $(\alpha, \theta)$ partition. Then the result is a regenerative $(\alpha, \theta)$ partition.*

PROOF. It is clear by construction that the split of table sizes into those to the left of table 1, table 1, and those to the right of table 1 is a Dirichlet$(\alpha, 1 - \alpha, \theta)$ split (cf. Lemma 5), and that given this split the dynamics of the composition to the left of table 1 and the composition to the right of table 1 produce limits as indicated. The conclusion now follows from the proposition. $\quad\square$

The particular cases $\theta = \alpha$ and $\theta = 0$ of Corollary 8 are known [23], Proposition 15. If $\theta = 0$, then $(G_1, D_1) = (G_1, 1)$ is the last component interval of $[0, 1] \setminus \mathcal{Z}_{\alpha,0}$ where $\mathcal{Z}_{\alpha,0}$ can be constructed as the restriction to $[0, 1]$ of the closed range of a stable subordinator of index $\alpha$. It is well known that the distribution of $G_1$ is then beta$(1 - \alpha, \alpha)$, and that the restriction of $\mathcal{Z}_{\alpha,0}$ to $[0, G_1]$ is a scaled copy of $\mathcal{Z}_{\alpha,\alpha}$ which can be defined by conditioning $\mathcal{Z}_{\alpha,0}$ on $1 \in \mathcal{Z}_{\alpha,0}$. Otherwise put, $\mathcal{Z}_{\alpha,0}$ and $\mathcal{Z}_{\alpha,\alpha}$ can be constructed as the zero sets of a Bessel process and standard Bessel bridge of dimension $2 - 2\alpha$. In the bridge case, $(G_1, D_1)$ can be represented as the interval covering a uniform random point independent of $\mathcal{Z}_{\alpha,\alpha}$, and $(G_1, D_1)$ splits $\mathcal{Z}_{\alpha,\alpha}$ into rescalings to $[0, G_1]$ and $[D_1, 1]$ of two independent copies of itself.

As indicated above, the local time process $(L(u), 0 \le u \le 1)$ can be described directly in terms of $\xi$ or $\mathcal{Z}_{\alpha,\theta}$: in the setting of Proposition 6, we have

$$(13) \qquad L(1 - \exp(-\xi_t)) = \Gamma(1 - \alpha) \int_0^t \exp(-\alpha \xi_s) \, ds;$$

cf. [14], Section 5. The right-continuous inverse of $L$ satisfies

$$(14) \quad L^{-1}(\ell) = 1 - \exp(-\xi_{T_\ell}) \qquad \text{where } T_\ell = \Gamma(1 - \alpha) \int_0^\ell \frac{dh}{(1 - L^{-1}(h))^\alpha}.$$

In fact, $(1 - L^{-1}(\ell), 0 \le \ell \le L(1))$ is a self-similar Markov process killed when reaching zero, so (13) and (14) are Lamperti's formulas [20] relating self-similar Markov processes and Lévy processes. This observation will tie in nicely with well-known properties of self-similar fragmentations that we introduce in Section 4.1. Furthermore, we will use the Stieltjes measure $dL^{-1}$ as a discrete measure on $[0, L(1)]$ to turn this interval into a string of beads in the sense of Definition 4.

2.4. *Finding the first table in the composition of table sizes.* Let $(\widetilde{\Pi}_n)$ be the sequence of random ordered partitions of $n$ induced by the ordered CRP, and $\mathcal{C}_n$ the regenerative composition structure of block sizes of $\widetilde{\Pi}_n$ studied in Proposition 6. According to (2), for each particular composition $(n_1, \ldots, n_\ell)$ of $n$,

$$(15) \qquad \mathbb{P}(\mathcal{C}_n = (n_1, \ldots, n_\ell)) = p_{\alpha,\theta}(n_1, \ldots, n_\ell)$$

$$:= \prod_{j=1}^\ell q_{\alpha,\theta}(N_j, n_j) \qquad \text{with } N_j := \sum_{i=j}^\ell n_i$$

for $q_{\alpha,\theta}$ as in (8). Now, for each $1 \le k \le \ell$, we wish to describe the conditional probability given this event that the first customer sits at the $k$th of these tables, which has size $n_k$.

LEMMA 9. *In the random ordered partition $\widetilde{\Pi}_n$ of $[n]$, given that the left-most block in this ordered partition is of size $n_1$, the probability that it contains 1 is*

$$(16) \qquad \frac{n_1 \theta}{n_1 \theta + N_2 \alpha} \qquad (N_2 := n - n_1).$$

*Given that the composition $\mathcal{C}_n$ of block sizes of $\widetilde{\Pi}_n$ is $(n_1, \ldots, n_\ell)$, for $1 \le k \le \ell$ the conditional probability that 1 falls in the $k$th block of size $n_k$ is*

$$(17) \qquad p_k^{(n)} \prod_{j=1}^{k-1} (1 - p_j^{(n)}) \qquad \text{for } p_j^{(n)} = \frac{n_j \theta}{n_j \theta + N_{j+1} \alpha} \text{ with } N_{j+1} := \sum_{i=j+1}^\ell n_i.$$

*In particular, if $\theta = \alpha$, then $\widetilde{\Pi}_n$ is exchangeable, and the size of the block containing 1 is a size-biased pick from the composition $\mathcal{C}_n$ of block sizes.*

PROOF. It is enough to describe the conditional probability, given that the first block has size $n_1$, that this block contains 1. For given that this block does not contain 1, the dynamics of the ordered CRP are such that the remainder of the ordered partition $\widetilde{\Pi}_n$, after order-preserving bijective relabeling (keeping label 1 fixed), makes a copy of $\widetilde{\Pi}_{n-n_1}$. The probability that the first block has size $n_1$ is found from (8) to be

$$(18) \qquad q_{\alpha,\theta}(n, n_1) = \binom{n-1}{N_2} \frac{[1-\alpha]_{n_1-1}}{[\theta + N_2]_{n_1}} \frac{(n_1 \theta + N_2 \alpha)}{n_1} \qquad (N_2 := n - n_1)$$

for $1 \leq n_1 \leq n$. In particular, for $n_1 = n$, the probability that there is just one block, $[n]$, is $[1-\alpha]_{n-1}/[1+\theta]_{n-1}$. This can also be seen directly from the sequential construction of the Chinese Restaurant. The denominator is the product of all weights for $n-1$ choices, and the numerator is the product of weights for each new customer sitting at the same table as all previous ones. The same direct argument shows that the probability that 1 ends up in the left-most block along with $n_1 - 1$ other integers is

$$(19) \qquad \binom{n-1}{N_2} \frac{[1-\alpha]_{n_1-1}[\theta]_{N_2}}{[1+\theta]_{n-1}},$$

where the first factor is the number of ways to choose which of the $n-1$ integers besides 1 are not in the first block, and, whatever this choice, the factors $[1-\alpha]_{n_1-1}$ and $[\theta]_{N_2}$ provide the product of weights of relevant remaining choices, and the denominator is the product of total weights. Look at the ratio of (19) and (18) to conclude. □

The case $\theta = 0$ deserves special mention. The probability of creating a new table to the right of the first $k$ tables is always zero. The effect of this is that 1 always remains in the right-most block of the ordered partition. Formula (16) in this case must be interpreted by continuity at $\theta = 0$, to give 0 for $1 \leq n_1 \leq n - 1$ and 1 for $n_1 = n$. This case is exceptional in that the size of the right-most table of the ordered restaurant has a strictly positive limiting proportion of all customers, with beta$(1 - \alpha, \alpha)$ distribution. This can be read, for example, from (4).

In all other cases the proportion at the right-most table converges almost surely to zero, as a consequence of (3). If $\alpha > 0$, the fraction in the left-most table tends to 0. If $\alpha = 0$ and $\theta > 0$, the fraction in the left-most table has a limiting beta$(1, \theta)$ distribution.

As $n$ tends to infinity, the rescaled compositions $\mathcal{C}_n$ become a limiting interval partition $\mathcal{Z}_{\alpha,\theta}$. Let us now study which interval of $\mathcal{Z}_{\alpha,\theta}$ is the limit of the block containing 1.

PROPOSITION 10. *Let $\theta > 0$. Given $\mathcal{Z}_{\alpha,\theta} = \{1 - \exp(-\xi_t), t \geq 0\}^{\mathrm{cl}}$, the conditional probability for the interval $(1 - \exp(-\xi_{t-}), 1 - \exp(-\xi_t))$ to be the limit of the block containing 1 is*

$$p(e^{-\Delta\xi_t}) \prod_{s<t}(1 - p(e^{-\Delta\xi_s})) \qquad \text{with } p(x) = \frac{(1-x)\theta}{(1-x)\theta + x\alpha}$$

*for all $t \geq 0$ with $\Delta\xi_t := \xi_t - \xi_{t-} > 0$.*

PROOF. For the random ordered partition $\widetilde{\Pi}_n = (\widetilde{\Pi}_n(1), \ldots, \widetilde{\Pi}_n(L_n(1)))$ and $u \in (0,1)$, we deduce from Lemma 9, in the notation of Proposition 6, that

$$\mathbb{P}(1 \in \widetilde{\Pi}_n(L_n(u))|\mathcal{Z}_{\alpha,\theta}^n) = p_{L_n(u)}^{(n)} \prod_{j=1}^{L_n(u)-1}(1 - p_j^{(n)}) \qquad \text{a.s.,}$$

where $\mathcal{Z}_{\alpha,\theta}^n := \{S_{n,j}/n, j \geq 0\} \to \mathcal{Z}_{\alpha,\theta} = [0,1] \setminus \bigcup_{i \in \mathcal{I}}(g_i, d_i)$ almost surely, with respect to the Hausdorff metric on closed subsets of $[0,1]$. We will refer to intervals $I_i = (g_i, d_i)$ as *parts* of $\mathcal{Z}_{\alpha,\theta}$. Denote $g_n(v) = \sup\{w \leq v : w \in \mathcal{Z}_{\alpha,\theta}^n\}$ and $d_n(v) = \inf\{w > v : w \in \mathcal{Z}_{\alpha,\theta}^n\}$ for $v \in (0,1)$, similarly, $g(v)$ and $d(v)$ for $\mathcal{Z}_{\alpha,\theta}$. For each fixed $v \in (0,1)$, we have

$$p_{L_n(v)}^{(n)} = \frac{(d_n(v) - g_n(v))\theta}{(d_n(v) - g_n(v))\theta + (1 - d_n(v))\alpha}$$

$$\to \frac{(d(v) - g(v))\theta}{(d(v) - g(v))\theta + (1 - d(v))\alpha} =: p_{g(v)} \qquad \text{a.s.}$$

Now fix $\varepsilon > 0$, then there is $M$ so that there are ("big") parts $I_1, \ldots, I_M$ of $\mathcal{Z}_{\alpha,\theta}$ that leave less than $\theta\varepsilon/8R$ uncovered, where $R = (1 - d(u))\alpha$. Using the a.s. convergence of left and right end points, a standard argument now shows that there is $N_0 \geq 0$ such that, for all $n \geq N_0$,

$$\left| \log(p_{L_n(u)}^{(n)}) + \sum_{j=1}^{L_n(u)-1} \log(1 - p_j^{(n)}) \right.$$

$$\left. - \log(p_{g(u)}) - \sum_{i \in \mathcal{I}: g_i < g(u)} \log(1 - p_{g_i}) \right| < \varepsilon,$$

since

$$\left| \log\left( \frac{(1-d(v))\alpha}{(d(v) - g(v))\theta + (1 - d(v))\alpha} \right) \right| \leq \frac{(d(v) - g(v))\theta}{(1 - d(v))\alpha} \leq \frac{(d(v) - g(v))\theta}{(1 - d(u))\alpha}$$

allows to jointly bound the sums of all small parts. Therefore,

$$\mathbb{P}(1 \in \widetilde{\Pi}_n(L_n(u))|\mathcal{Z}_{\alpha,\theta}^n) = p_{L_n(u)}^{(n)} \prod_{j=1}^{L_n(u)-1}(1 - p_j^{(n)})$$

$$\to p_{g(u)} \prod_{i \in \mathcal{I}: g_i < g(u)} (1 - p_{g_i}) \qquad \text{a.s.}$$

Now we use dominated convergence to deduce for any bounded continuous function $f$ on the space of closed subsets of $[0, 1]$ (equipped with the Hausdorff metric) that

$$\mathbb{E}(f(\mathcal{Z}_{\alpha,\theta}^n)\mathbb{P}(1 \in \widetilde{\Pi}_n(L_n(u)) | \mathcal{Z}_{\alpha,\theta}^n)) \to \mathbb{E}\left(f(\mathcal{Z}_{\alpha,\theta}) p_{g(u)} \prod_{i \in \mathcal{I}: g_i < g(u)} (1 - p_{g_i})\right).$$

However, we also have

$$\mathbb{E}(f(\mathcal{Z}_{\alpha,\theta}^n) 1_{\{1 \in \widetilde{\Pi}_n(L_n(u))\}}) \to \mathbb{E}(f(\mathcal{Z}_{\alpha,\theta}) 1_{\{u \in (G_1, D_1)\}})$$
$$= \mathbb{E}(f(\mathcal{Z}_{\alpha,\theta})\mathbb{P}(u \in (G_1, D_1) | \mathcal{Z}_{\alpha,\theta})),$$

since the distributions of $G_1$ and $D_1$ are continuous or degenerate ($G_1 = 0$ or $D_1 = 1$) by Corollary 8. We identify

$$(20) \qquad p_{g(u)} \prod_{i \in \mathcal{I}: g_i < g(u)} (1 - p_{g_i})$$

as a version of the conditional probability $\mathbb{P}(u \in (G_1, D_1) | \mathcal{Z}_{\alpha,\theta})$ for all $u \in (0, 1)$.

Finally, conditionally given $\mathcal{Z}_{\alpha,\theta}$, each of the *countable* number of times $t$ such that $\xi_{t-} < \xi_t$ is associated with an interval $(1 - \exp(-\xi_{t-}), 1 - \exp(-\xi_t))$ of $u$-values to which (20) applies, so the conditional distribution of $(G_1, D_1)$ given $\mathcal{Z}_{\alpha,\theta}$ is as claimed. $\quad\square$

The limiting interval in $\mathcal{Z}_{\alpha,\theta}$ of the block containing 1 corresponds to a jump of the (for $\theta = 0$ killed by an infinite jump at an exponential time **e**) subordinator $\xi$. Denote the time of this jump by $\tau$. It can now be checked directly that the boundary points $(1 - \exp(-\xi_{\tau-}), 1 - \exp(-\xi_\tau))$ describe a Dirichlet$(\alpha, 1 - \alpha, \theta)$ split of $[0, 1]$ as shown in Corollary 8. Standard thinning arguments for the Poisson point process $(\Delta\xi_t, t \geq 0)$ show that $\xi_{\tau-} \stackrel{d}{=} \xi_\tau^0$, where $\xi^0$ is a subordinator independent of $\tau$ with Lévy measure $(1 - p(e^{-x}))\Lambda_{\alpha,\theta}(dx)$ and Laplace exponent

$$\Phi_0(s) = \int_0^\infty (1 - e^{-sx})(1 - p(e^{-x}))\Lambda_{\alpha,\theta}(dx)$$

so that

$$\mathbb{E}(e^{-s\xi_{\tau-}}) = \int_0^\infty e^{-t\Phi_0(s)} \lambda e^{-\lambda t} \, dt = \frac{\lambda}{\Phi_0(s) + \lambda},$$

where $\lambda = \Gamma(1 - \alpha)\Gamma(\theta + 1)/\Gamma(\theta + 1 - \alpha)$ is the rate of the exponential variable $\tau$.

For simplicity, let $\theta > 0$. The case $\theta = 0$ is similar, taking into account the killing at the infinite jump. We find the Lévy measure $\Lambda_{\alpha,\theta}(dx)$ of $\xi$ from $\Phi_{\alpha,\theta}(s) = \int_{(0,\infty)}(1 - e^{-sx})\Lambda_{\alpha,\theta}(dx)$ with $\Phi_{\alpha,\theta}$ given in (10) (cf. also [13], formula (41)) and change variables $u = e^{-x}$ to get

$$\Phi_0(s) = \Phi_{\alpha,\theta}(s) - \theta \int_0^1 (1 - u^s) u^{\theta-1}(1 - u)^{-\alpha}\, du$$

$$= sB(s + \theta, 1 - \alpha) - \theta(B(\theta, 1 - \alpha) - B(s + \theta, 1 - \alpha))$$

$$= (s + \theta)B(\theta + s, 1 - \alpha) - \lambda, \qquad \text{where } B(a, b) = \Gamma(a)\Gamma(b)/\Gamma(a + b),$$

and, hence,

$$\mathbb{E}(e^{-s\xi_{\tau-}}) = \frac{\theta}{\theta + s} \frac{\Gamma(\theta)\Gamma(s + \theta + 1 - \alpha)}{\Gamma(\theta + 1 - \alpha)\Gamma(s + \theta)}.$$

These are the moments of a beta$(\alpha, \theta + 1 - \alpha)$ distribution in accordance with Corollary 8. Similarly, $\Delta\xi_\tau$ has distribution

$$\frac{1}{\lambda}p(e^{-x})\Lambda_{\alpha,\theta}(dx)$$

and so the interval size $\exp(-\xi_{\tau-})(1 - \exp(-\Delta\xi_\tau))$ relative to the remaining proportion $\exp(-\xi_{\tau-})$ can be seen to be independent of $\exp(-\xi_{\tau-})$ and to have a beta$(1 - \alpha, \theta)$ distribution. By Lemma 5(b), this establishes the Dirichlet$(\alpha, 1 - \alpha, \theta)$ distribution of Corollary 8.

## 3. Markov branching models and weighted discrete $\mathbb{R}$-trees with edge lengths.

### 3.1. *Markov branching models.*
Our formalism for combinatorial trees follows [18], Section 2. For $n = 1, 2, \ldots,$ let $T_n^\circ$ denote a random unlabeled rooted binary tree with $n$ leaves. The sequence $(T_n^\circ, n \geq 1)$ is said to have the *Markov branching property* [2, 10] if conditionally given that the first split of $T_n^\circ$ is into tree components whose numbers of leaves are $m$ and $n - m$, these components are like independent copies of $T_m^\circ$ and $T_{n-m}^\circ$, respectively. The distributions of the first splits of $T_n^\circ$, $n \geq 1$, are denoted by $(q^\circ(m, n - m), 1 \leq m \leq n/2)$ and referred to as the splitting rule of $(T_n^\circ, n \geq 1)$.

For a finite set $B$, let $\mathbb{T}_B$ be the set of binary trees with leaves labeled by $B$. For $T_n \in \mathbb{T}_{[n]}$ and $B \subset [n]$, let $T_{n,B} \in \mathbb{T}_B$ be the reduced subtree of $T_n$ spanned by leaves in $B$, and let $\widetilde{T}_{n,B} \in \mathbb{T}_{[\#B]}$ be the image of $T_{n,B}$ after relabeling of leaves by the increasing bijection from $B$ to $[\#B]$. It will be convenient to label each branch point of $T_n$ by the set of leaf labels in the subtree above the branch point. A tree $T_n \in \mathbb{T}_{[n]}$ is then uniquely represented by a collection of subsets of $[n]$. Such a tree has the natural interpretation as a fragmentation tree, where blocks (i.e. labels of branch points, $[n]$ for the first branch point) fragment as one passes from one level to the next. We will write $B \in T_n$ if $T_n$ has a vertex with label $B$.

PROPOSITION 11. *Let $(T_n^{\alpha,\theta}, n \geq 1)$ for some $0 \leq \alpha \leq 1$ and $\theta \geq 0$ be an $(\alpha,\theta)$-tree growth process as defined in Definition 3. Then:*

(a) *the delabeled process $(T_n^{\alpha,\theta,\circ}, n \geq 1)$ has the Markov branching property with splitting rule*

$$q^\circ(m, n-m) = q_{\alpha,\theta}(n-1, m) + q_{\alpha,\theta}(n-1, n-m), \qquad 1 \leq m < n/2,$$

$$q^\circ(n/2, n/2) = q_{\alpha,\theta}(n-1, n/2), \qquad \textit{if } n \textit{ is even},$$

*where $q_{\alpha,\theta}(n, m)$ is given in (8);*

(b) *the labeled process $(T_n^{\alpha,\theta}, n \geq 1)$ is regenerative in the sense that for each $n \geq 1$, conditionally given that the first split of $T_n^{\alpha,\theta}$ is by a partition $\{B, [n] \setminus B\}$ of $[n]$ with $\#B = m$, the relabeled subtrees $\tilde{T}_{n,B}^{\alpha,\theta}$ and $\tilde{T}_{n,[n] \setminus B}^{\alpha,\theta}$ are independent copies of $T_m^{\alpha,\theta}$ and $T_{n-m}^{\alpha,\theta}$, respectively.*

PROOF. For notational convenience, we drop superscripts $\alpha, \theta$. Recall from the Introduction the identification (1) of leaf $k+1$ of $(T_n, n \geq 1)$ and customer $k$ of the regenerative composition structure $(\mathcal{C}_n, n \geq 1)$ of the ordered Chinese Restaurant Process described in Proposition 6, for all $k \geq 1$. This identifies $\mathcal{C}_{n-1}$ as the composition of subtree sizes growing off the spine from the root to leaf 1. In particular, we see that for each $n \geq 2$ the distribution $q^\circ$ stated here applies as splitting rule at the first branch point of $T_n$ and indeed on the spine of $T_n$.

To establish the Markov branching property, proceed by induction. $T_1^\circ$, $T_2^\circ$ and $T_3^\circ$ trivially have the Markov branching property. Assume that the property is true for $T_1^\circ, \ldots, T_n^\circ$ for some $n \geq 3$. Then, by the growth procedure, two scenarios can occur. Given $n+1$ attaches to the trunk, the subtrees of $T_{n+1}^\circ$ are $T_n^\circ$ and the deterministic tree with single leaf $n+1$, they are trivially conditionally independent and, by the induction hypothesis, have distributions as required. Given $n+1$ attaches in one or the other subtree of $T_n^\circ$ of sizes $m$ and $n-m$, the induction hypothesis yields the conditional independence and Markov branching distributions for these subtrees, and also yields that the insertion of a new leaf into one of these trees gives the corresponding Markov branching distribution of size $m+1$ or $n-m+1$, respectively, by the recursive nature of the growth procedure.

This proves (a). The induction is easily adapted to also prove (b). Just note that the $(\alpha, \theta)$-tree growth rules are invariant under increasing bijections from $B$ to $[\#B]$. □

3.2. *Sampling consistency and the proof of Proposition 1.* Recall that a sequence of trees $(T_n^\circ, n \geq 1)$ is *weakly sampling consistent* if uniform random removal of a leaf of $T_{n+1}^\circ$ yields a reduced tree with the same distribution as $T_n^\circ$, for each $n \geq 1$.

For $(T_n^{\alpha,\theta}, n \geq 1)$ with splitting rules $q^\circ(m, n-m)$ as before (with $m \leq n-m$), to match notation with Ford [10], Proposition 41, introduce the split probability functions

- $q^{\text{bias}}(x,y)$ defined so that $q^{\text{bias}}(m, n-m) = q_{\alpha,\theta}(n-1, m)$ [see (8)] is the probability that $[n]$ is first split into pieces of size $m$ and $n-m$, for $1 \leq m \leq n-1$, where we are supposing that the piece of size $m$ does not contain label 1; so $q^{\text{bias}}(x,y) = q_{\alpha,\theta}(x + y - 1, x)$;
- $q^{\text{sym}}(x,y) = \frac{1}{2}q^{\text{bias}}(x,y) + \frac{1}{2}q^{\text{bias}}(y,x)$ for the symmetrization of $q^{\text{bias}}$. Then we have $q^{\text{sym}}(x,y) = \frac{1}{2}q^\circ(x,y)$ for all $x < y$ and $q^{\text{sym}}(x,x) = q^\circ(x,x) = q^{\text{bias}}(x,x)$ for all $x \geq 1$.

Ford uses symmetrized splitting rules to grow unlabeled planar trees. For us they are useful for a weak sampling consistency criterion: let

$$d^{\text{sym}}(x,y) := q^{\text{sym}}(x,y)\left(1 - \frac{q^{\text{sym}}(1, x+y) + q^{\text{sym}}(x+y, 1)}{x+y+1}\right)$$
$$- q^{\text{sym}}(x+1, y)\frac{x+1}{x+y+1} - q^{\text{sym}}(x, y+1)\frac{y+1}{x+y+1}.$$

Ford [10], Proposition 41, showed that $(T_n^\circ)$ is weakly sampling consistent if and only if $d^{\text{sym}}(x,y) = 0$ for all positive integers $x$ and $y$. He verified this property for the $(\alpha, 1-\alpha)$-trees.

PROOF OF PROPOSITION 1(c). For the $(\alpha, \theta)$ splitting rules we obtain

$$d^{\text{sym}}(1,1) = d^{\text{sym}}(1,2) = 0,$$

but

$$(21) \quad d^{\text{sym}}(1,3) = \frac{(1-\alpha)(1-\alpha-\theta)(2-\alpha-\theta)(3-\alpha+\theta)(\alpha+\theta)}{10(1+\theta)^2(2+\theta)^2(3+\theta)},$$

which shows that a necessary condition for $(T_n^\circ)$ to be weakly sampling consistent is that $\theta$ equals either $1 - \alpha$ or $2 - \alpha$. Ford showed that $\theta = 1 - \alpha$ produces weakly sampling consistent trees. The proof of part (c) of Proposition 1 is completed by the following lemma. □

LEMMA 12. *For $\theta = 1 - \alpha$ and $\theta = 2 - \alpha$, the symmetrized splitting rules are the same. Therefore, the $(\alpha, 2-\alpha)$ tree growth process is weakly sampling consistent.*

PROOF. For convenience of notation in this proof, denote the nonsymmetric splitting rules, for Ford's case $\theta = 1 - \alpha$ by

$$q_n^F(m, n-m) = \binom{n-1}{m} \frac{m + (n-1-2m)\alpha}{n-1} \frac{\Gamma(m-\alpha)\Gamma(n-m-\alpha)}{\Gamma(1-\alpha)\Gamma(n-\alpha)}$$

[see (8)], and for $\theta = 2 - \alpha$ by

$$q^X(m, n-m) = \binom{n-1}{m} \frac{2m + (n-1-2m)\alpha}{n-1} \frac{\Gamma(m-\alpha)\Gamma(n-m+1-\alpha)}{\Gamma(1-\alpha)\Gamma(n+1-\alpha)}.$$

Now the claim is that

$$\tfrac{1}{2}q^X(m, n-m) + \tfrac{1}{2}q^X(n-m, m) = \tfrac{1}{2}q^F(m, n-m) + \tfrac{1}{2}q^F(n-m, m),$$

which after the obvious cancellations is equivalent to

$$(n-m)(2m + (n-1-2m)\alpha)(n-m-\alpha)$$
$$+ m(2n - 2m + (2m-n-1)\alpha)(m-\alpha)$$
$$= (n-m)(m + (n-1-2m)\alpha)(n-\alpha)$$
$$+ m(n-m + (2m-n-1)\alpha)(n-\alpha),$$

and this is easily checked.  □

The nonsymmetrized rules are equal only if $\alpha = 1$, trivially, since this is the deterministic comb model, where all leaves connect to a single spine. In fact, it can be shown that these coincidences of symmetrized splitting rules are the only such coincidences, in particular, for fixed $\alpha$, the splitting rules as a path in the space of splitting rules, parameterized by $\theta \geq 0$, have precisely one loop.

Let us turn to strong sampling consistency and exchangeability.

PROOF OF PROPOSITION 1(a)–(b). Assume that $(T_n, n \geq 1)$ is strongly sampling consistent for some $\theta \in \{1 - \alpha, 2 - \alpha\}$, then it is not hard to show that also the regenerative composition structure $(\mathcal{C}_n, n \geq 1)$ generated by the associated ordered Chinese Restaurant Process is strongly sampling consistent. By Proposition 6, this implies $\theta = \alpha$ and, hence, $\theta = \alpha = 1/2$. On the other hand, it is well known that this case is strongly sampling consistent. This establishes part (b) of Proposition 1.

Part (a) of Proposition 1 is easily checked for $n = 3$. The shape $T_3^\circ$ is deterministic, as there is only one rooted binary tree with three leaves. This tree has one leaf at height 2 and two leaves at height 3. Denote the label of the leaf at height 2 by $M$. Then exchangeability requires

$$\frac{1}{3} = \mathbb{P}(M = 2) = \frac{\theta}{1 + \theta} \quad \Rightarrow \quad \theta = \frac{1}{2}$$

and for $\theta = 1/2$,

$$\frac{1}{3} = \mathbb{P}(M = 3) = \frac{\alpha}{1 + \theta} = \frac{2\alpha}{3} \quad \Rightarrow \quad \alpha = \frac{1}{2},$$

using the growth rules. This completes the proof of Proposition 1. $\square$

We conclude this subsection by a study of boundary cases. For $\alpha = 1$, we have a comb model (all leaves directly attached to a single spine) with nonuniform labeling (for $\theta = 1$, leaves $2, 3, \ldots$ are exchangeable, and for $\theta = 0$, leaves $3, 4, \ldots$ are exchangeable), but strongly sampling consistent as the delabeled trees are deterministic. The trees grow linearly in height.

For $\alpha = 0$, we get a tree growth model that one might call internal boundary aggregation on the complete binary tree in a beta$(1, \theta)$ random environment. Informally, attach $n + 1$ to $T_n$ at the terminal state of a walker climbing the tree by flipping the beta$(1, \theta)$ coin corresponding to each branch point until he reaches a leaf of $T_n$. Insert $n + 1$ by replacing the leaf by a new branch point connected to the leaf and $n + 1$.

More formally, let $\mathcal{X} = \bigcup_{n \geq 0} \{0, 1\}^n$ be the complete rooted binary tree, where $\{0, 1\}^0 = \varnothing$ is the empty word and elements of $\{0, 1\}^n$ are identified as binary words of length $n$. Mark all vertices of $\mathcal{X}$ by independent beta$(1, \theta)$ random variables $W_x$, $x \in \mathcal{X}$. Consider the binary tree growth process with edge selection rule as follows:

(i)$^{\mathcal{W}}$ Let a walker start from $Z_0 = [n]$, with $X_0 = \varnothing$ (for $k = 0$), with steps as in (ii)$^{\mathcal{W}}$.

(ii)$^{\mathcal{W}}$ Given $T_n$ and a word $X_k$, let $X_{k+1} \sim \text{Bernoulli}(W_{X_k})$. If $X_{k+1} = 1$ and $Z_k$ has children $B$ and $Z_k \setminus B$, where $B$ contains the smallest label of $Z_k$, set $Z_{k+1} = B$, otherwise $Z_{k+1} = Z_k \setminus B$. If $\#Z_{k+1} \geq 2$, repeat (ii)$^{\mathcal{W}}$. Otherwise select edge $Z_{k+1} = \{L_{n+1}\}$.

In our formalism where $T_n$ is a collection of subsets of $[n]$, the growth step can be made explicit as $T_{n+1} = \{B \cup \{n+1\} : L_{n+1} \in B \in T_n\} \cup \{B : L_{n+1} \notin B \in T_n\} \cup \{\{L_{n+1}\}, \{n+1\}\}$.

PROPOSITION 13. (a) *The family* $(T_n)_{n \geq 1}$ *grown via* (i)$^{\mathcal{W}}$–(ii)$^{\mathcal{W}}$ *is a* $(0, \theta)$-*tree growth process.*

(b) *The labeling of* $T_n$, $n \geq 3$ *is not exchangeable for any* $\theta \geq 0$; *the trees are weakly sampling consistent if and only if* $\theta = 0$ *or* $\theta = 1$ *or* $\theta = 2$; *the trees grow logarithmically (except for* $\theta = 0$, *when the model is the comb model and growth is linear).*

PROOF. (a) This follows directly from the growth rules of the $(0, \theta)$-tree growth process, since internal edges are never selected for insertions. The first branch point separates 1 and 2. At this branch point, and inductively every other branch point, an urn scheme governs the selection procedure, with initial weight 1 for the subtree of the larger label, $\theta$ for the subtree of the smaller label, so a beta$(1, \theta)$ limiting proportion of insertions will take place in the subtree of the larger label; cf. Lemma 5.

(b) The exchangeability claim follows easily from the growth procedures. Weak sampling consistency can be read from (21), which also holds for $\alpha = 0$. Logarithmic growth follows from the following considerations. Just as we argued for $0 < \alpha < 1$ in the Introduction, also for $\alpha = 0$, the height $K_n$ of leaf 1 in $T_n$ has the same dynamics as the number of tables in a Chinese Restaurant Process with $(0, \theta)$ seating plan. In this case $K_n$ is known to grow logarithmically, with $K_n / \log(n) \to \theta$ if $\theta > 0$. It is easy to see that also the rescaled height of leaf $k$ converges to $\theta$.  □

Note that the height of the branch point between any two leaves $j$ and $k$ is constant, hence converges to zero when rescaled by $\log(n)$. Therefore, in a logarithmically scaled limit tree all leaves would be adjacent to the root with no further branching structure.

3.3. *Weighted discrete $\mathbb{R}$-trees with edge lengths.* A pointed compact metric space $(\mathcal{T}, d, \rho)$ is called a *compact $\mathbb{R}$-tree* with root $\rho \in \mathcal{T}$ if it is complete separable path-connected and has the tree property:

- for any $\sigma, \sigma' \in \mathcal{T}$ there is a unique isometry $g_{\sigma,\sigma'} : [0, d(\sigma, \sigma')] \to \mathcal{T}$ such that $g_{\sigma,\sigma'}(0) = \sigma$ and $g_{\sigma,\sigma'}(d(\sigma, \sigma')) = \sigma'$; denote $[[\sigma, \sigma']] = g_{\sigma,\sigma'}([0, d(\sigma, \sigma')])$; furthermore, any simple path from $\sigma$ to $\sigma'$ has range $[[\sigma, \sigma']]$.

In this section we restrict our attention to $\mathbb{R}$-tree representatives of discrete trees with edge lengths such as $T_n \in \mathbb{T}_{[n]}$ with edge lengths $e_B \in (0, \infty)$, $B \in T_n$, where $e_B$ refers to the parent edge *below* $B$, so $e_{[n]}$ is the length of the root edge. For $B \in T_n$ with ancestors $[n] = B_0 \supset B_1 \supset \cdots \supset B_k = B$ in $T_n$, we denote its birth time by $l_B = e_{B_0} + \cdots + e_{B_{k-1}}$ and its death time by $r_B = l_B + e_B$. Recall, for example, from [7] that we can associate a real tree

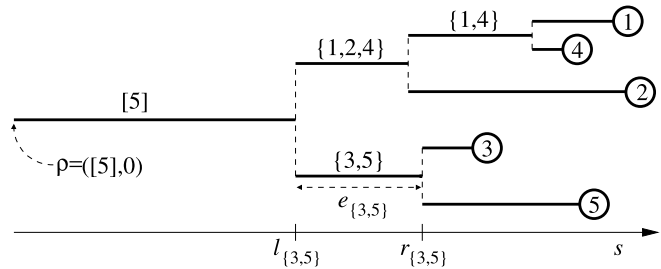

Fig. 4.   *Canonical representation of a tree $T_5$ with edge lengths $e_B$, $B \in T_5$.*

as a subset of $T_n \times [0, \infty)$ as

$$(22) \qquad \mathcal{T} = \{([n], 0)\} \cup \{(B, s) : B \in T_n, s \in (l_B, r_B]\},$$

in *canonical form*, so that $E_B := B \times (l_B, r_B]$ represents the edge below $B$ of *Euclidean* length $e_B = r_B - l_B$; cf. Figure 4. We refer to $T_n$ as the *shape* of $\mathcal{T}$. We define the root $\rho = ([n], 0)$ and a metric $d$ on $\mathcal{T}$ that extends the natural Euclidean metric on the edges and that connects the edges to a tree. If $\sigma = (B, s) \in \mathcal{T}$, then we set $d(\rho, \sigma) = s$. Let $\sigma' = (B', s') \in \mathcal{T} \setminus \{\rho\}$. We define $d(\sigma, \sigma')$ by

$$d(\sigma, \sigma') = \begin{cases} d(\rho, \sigma) + d(\rho, \sigma') - 2r_{B \vee B'}, \\ \qquad \text{if } B \vee B' := \bigcap_{B'' \in T_n : B \cup B' \subset B''} B'' \notin \{B, B'\}; \\ |d(\rho, \sigma) - d(\rho, \sigma')|, \\ \qquad \text{otherwise, that is, if } B \subseteq B' \text{ or } B' \subseteq B; \end{cases}$$

here the first case is when $B \cap B' = \varnothing$, that is, there is a branch point, the *last common ancestor* $B \vee B'$, for which $B$ is in one subtree and $B'$ in the other.

A *weighted* $\mathbb{R}$-*tree* is equipped with a probability measure $\mu$ on the Borel sets of $(\mathcal{T}, d)$. As a relevant example consider an interval partition $\mathcal{Z} \subset [0, 1]$ with local time $(L(u), 0 \le u \le 1)$. We can associate a real tree consisting of a single branch $[0, L(1)]$ and specify $\mu$ by its distribution function $L^{-1}$, that is, $\mu([0, L(u)]) = u$. We visualize the atoms of different sizes lined up on $[0, L(1)]$ (particularly if they are dense, but also if they are not dense) as a *string of beads* and use this term to refer to the weighted interval; cf. Figure 3 in the Introduction for a tree composed of strings of beads. In this specific single-branch context we have a natural notion of convergence, namely, weak convergence of Stieltjes measures $dL^{-1}$ as measures on $[0, \infty)$, where the interval $[0, L(1)]$ is determined by the supremum of the support of the measure. In this sense, Proposition 6 easily yields the following convergence of strings of beads:

$$(23) \qquad ([0, n^{-\alpha} L_n(1)], d(n^{-\alpha} L_n)^{-1}) \to ([0, L(1)], dL^{-1}) \qquad \text{weakly a.s.}$$

In general, $\mathbb{R}$-trees can have features such as a dense set of *branch points* ($\sigma \in \mathcal{T}$ such that $\mathcal{T} \setminus \{\sigma\}$ has three or more connected components) and allow diffuse weight measures on an uncountable set of *leaves* ($\sigma \in \mathcal{T}$ such that $\mathcal{T} \setminus \{\sigma\}$ is connected). We will introduce a suitable space of $\mathbb{R}$-trees and the weighted Gromov–Hausdorff notion of convergence in Section 4.1, self-similar fragmentation trees will be introduced as relevant examples.

3.4. *Convergence of reduced trees and the proof of Proposition 2.* Recall from the Introduction our notation $R(T_n; [k])$ for the reduced tree, the subtree of $T_n$ spanned by leaves labeled $[k]$ and equipped with the graph distances in $T_n$ as edge lengths. We now associate an $\mathbb{R}$-tree via (22). Proposition 2 claims that the $(\alpha, \theta)$-tree growth process $(T_n, n \geq 1)$ has the asymptotics

$$n^{-\alpha} R(T_n, [k]) \to \mathcal{R}_k \qquad \text{in the Gromov–Hausdorff sense as } n \to \infty$$

for some limiting discrete $\mathbb{R}$-tree $\mathcal{R}_k$ with random edge lengths and precisely $k$ leaves labeled by $[k]$. To describe the distribution recursively, we will use notation $\mathcal{S}_k^B = \{(B, l_B)\} \cup \{(A, s) \in \mathcal{R}_k : A \subset B\}$ for the *subtree* of $\mathcal{R}_k$ above $B$. In the following Proposition 14 we prove a refinement of Proposition 2 that includes a mass measure $\mu_k$ on the branches of $\mathcal{R}_k$.

DEFINITION 5. Let $(\mathcal{S}, d_k|_{\mathcal{S}}, \mu_k|_{\mathcal{S}})$ be a closed connected subset of $(\mathcal{R}_k, d_k, \mu_k)$ with mass $m = \mu_k(\mathcal{S}) > 0$ and root $(B, s_0)$ given by $B = \bigcup_{(A,s) \in \mathcal{S}} A$, $s_0 = \min\{s : (A, s) \in \mathcal{S}\}$. Then we associate the *relabeled, scaled and shifted* tree $(\widetilde{\mathcal{S}}, \widetilde{d}, \widetilde{\mu})$ as the canonical form (22) of the tree $\mathcal{S}$ with edge lengths multiplied by $m^{-\alpha}$, labels changed by the increasing bijection from $B$ to $[\#B]$, mass measure pushed forward via these operations and then multiplied by $m^{-1}$.

Once we have embedded $\mathcal{R}_k$ as a subtree of a CRT $(\mathcal{T}, \rho, \mu)$, the atoms of the mass measure $\mu_k$ will correspond to the $\mu$-masses of the connected components of $\mathcal{T} \setminus \mathcal{R}_k$ projected onto $\mathcal{R}_k$. More formally, for any two $\mathbb{R}$-trees $\mathcal{R} \subset \mathcal{T}$ with common root $\rho \in \mathcal{R}$, there is a natural projection

$$\pi^{\mathcal{R}} : \mathcal{T} \to \mathcal{R}, \qquad u \mapsto g_{\rho,\sigma}(\sup\{t \geq 0 : g_{\rho,\sigma}(t) \in \mathcal{R}\}),$$

where $g_{\rho,\sigma} : [0, d(\rho, \sigma)] \to \mathcal{T}$ is the unique isometry with $g_{\rho,\sigma}(0) = \rho$ and $g_{\rho,u}(d(\rho, \sigma)) = \sigma$. For a measure $\mu$ on $\mathcal{T}$, we denote the push-forward via $\pi^{\mathcal{R}}$ by

$$\pi_*^{\mathcal{R}} \mu(C) = (\pi^{\mathcal{R}})^{-1}(C), \qquad C \in \mathcal{B}(\mathcal{R}) := \{D \subset \mathcal{R} : D \text{ Borel measurable}\}.$$

Denote by $\nu_n$ the empirical (probability) measure on the leaves of the $\mathbb{R}$-tree representation of $T_n$ with unit edge lengths. We refer to $\nu_n$ as *mass measure* of $T_n$.

PROPOSITION 14. *Denote by $(T_n, n \geq 1)$ an $(\alpha, \theta)$-growth process as defined in Definition 3.*

(a) *Let $0 < \alpha < 1$, $\theta \geq 0$ and $k \geq 1$. We have, as $n \to \infty$, that*

$$(24) \qquad (n^{-\alpha} R(T_n, [k]), \pi_*^{R(T_n, [k])} \nu_n) \to (\mathcal{R}_k, \mu_k) \qquad \text{weakly a.s.}$$

*in the sense that for all $2k - 1$ edges the strings of beads converge a.s. as in (23).*

(b) *Let $0 < \alpha < 1$ and $\theta > 0$. The distribution of $(\mathcal{R}_k, \mu_k)$ is determined recursively as follows. $(\mathcal{R}_1, \mu_1) = (E_{\{1\}}, \mu_1)$ is an $(\alpha, \theta)$-string of beads. For $k \geq 2$, $(\mathcal{R}_k, \mu_k)$ has shape $T_k$ and the first branch point splits $(\mathcal{R}_k, \mu_k)$ into three components: a trunk and two subtrees. Conditionally given that $T_k$ first branches into $\{B, [k] \setminus B\}$ with $1 \in [k] \setminus B$ and $\#B = m$, the following four random variables are independent:*

- *$(H_1, H_2, H_3) = (\mu_k(E_{[k]}), \mu_k(\mathcal{S}_k^B), \mu_k(\mathcal{S}_k^{[k]\setminus B})) \sim \text{Dirichlet}(\alpha, m - \alpha, k - m - 1 + \theta)$;*
- *the scaled and shifted trunk $(\widetilde{E}_{[k]}, \widetilde{\mu}_k^{E_{[k]}})$ is an $(\alpha, \alpha)$-string of beads;*
- *the relabeled, scaled and shifted subtree $(\widetilde{\mathcal{S}}_k^B, \widetilde{\mu}_k^B)$ is distributed as $(\mathcal{R}_m, \mu_m)$,*
- *the relabeled, scaled and shifted subtree $(\widetilde{\mathcal{S}}_k^{[k]\setminus B}, \widetilde{\mu}_k^{[k]\setminus B})$ as $(\mathcal{R}_{k-m}, \mu_{k-m})$.*

PROOF. The proof is an extension of the proof of [18], Proposition 18. The case $k = 1$ was established in (23). Now fix $k \geq 2$ and $T_k$. Assume, inductively, that the proposition is proved up to tree size $k - 1$. For $n \geq k$, the reduced trees $(R(T_n, [k]), \pi_*^{R(T_n, [k])} \nu_n)$ all have the same shape as $T_k$. In the transition from $n$ to $n + 1$, mass increases by 1, and there may be no change of the reduced tree, or one of the edge lengths may increase by 1.

Let us first just distinguish the weights of the trunk below the first branch point and the two subtrees above, of sizes $m$ and $k - m$, say. We can associate three colors with the three components. It is easy to see that the mass allocation behaves like an urn model. The $(\alpha, \theta)$-tree growth rules specify initial urn weights of $\alpha$, $m - \alpha$ and $k - m - 1 + \theta$. Hence, these are the parameters of the Dirichlet distribution of limiting urn weights $(H_1, H_2, H_3)$; cf. Lemma 5.

Now we can treat separately the evolution of the three components, conditionally given $(H_1, H_2, H_3)$. See the proof of [18], Proposition 18, for details of this argument, which gives us the claimed independence.

The trunk follows the dynamics of an $(\alpha, \alpha)$ ordered CRP (when restricted to the proportion $H_1$ of leaves added in this part of the tree) whose limiting behavior was studied in Proposition 6 and (23). By the recursive nature of the growth procedure, the two subtrees have the same dynamics as $(R(T_n, [m]), \pi_*^{R(T_n, [m])} \nu_n)$ and $(R(T_n, [k-m]), \pi_*^{R(T_n, [k-m])} \nu_n)$, respectively, (when restricted to the proportions $H_2$ and $H_3$ of leaves added to these parts), and the induction hypothesis establishes their limiting behavior. □

PROOF OF PROPOSITION 2. Joint convergence with mass measures in Proposition 14(a) implies convergence of the trees without mass measures, so the proof of Proposition 2 is complete. □

The result in (b) is still true for $\theta = 0$, if interpreted appropriately. In fact, leaf edges with zero edge weight disappear in the limit of (a). It is now implicit in the above description that the limits of the associated leaves are on branches of the limiting tree. They are not leaves themselves. In particular, the first split is not necessarily at the first (topological) branch point of $(\mathcal{R}_k, \mu_k)$, but (for $m = k - 1$) may be leaf 1 on the branch leading to the first (topological) branch point. If so, it is this splitting the recursive description describes, with zero mass proportion for the degenerate subtree containing 1 (zero third parameter for the Dirichlet distribution).

3.5. *Growth of $(\mathcal{R}_k, \mu_k)$ by bead crushing.* The recursion can be partially solved to give the distribution of $(\mathcal{R}_k, \mu_k)$ more explicitly. Specifically, standard Dirichlet calculations [e.g., using Lemma 5(b)] show that the mass splits introduced by the branch points on the spine from the root to 1 lead to Dirichlet mass splits with parameter $\theta$ for the edge adjacent to 1, parameter $\alpha$ for all other spinal edges and parameter $m - \alpha$ for every subtree with $m$ leaves. When applying the recursion in a subtree off the spine with $m$ leaves, we have $m - \alpha = m - 1 + \theta$ only if $\theta = 1 - \alpha$, so only in the $(\alpha, 1 - \alpha)$ case, the overall mass split edge by edge is Dirichlet distributed, Dirichlet$(\alpha, \ldots, \alpha, 1 - \alpha, \ldots, 1 - \alpha)$ with $\alpha$ for the $n - 1$ inner branches and $1 - \alpha$ for the $n$ leaf edges. For $\theta \neq 1 - \alpha$, we get a mass split edge by edge that is best described recursively. Regarding the mass distribution on edges, we note:

COROLLARY 15. *In the setting of Proposition 14, conditionally given $T_k$ and an edge-by-edge split*

$$(\mu_k(E_B), B \in T_k) = (h_B, B \in T_k),$$

*the components $(E_B, \mu_k|_{E_B})$ are independent and such that $(\widetilde{E}_B, \widetilde{d}_k^{E_B}, \widetilde{\mu}_k^{E_B})$ is an $(\alpha, \alpha)$-string of beads for $\#B \geq 2$ and an $(\alpha, \theta)$-strings of beads for $\#B = 1$.*

Since the Dirichlet mass proportions induced by the split at the first branch point are independent from the three rescaled components in Proposition 14(b), the $(\alpha, \theta)$-tree growth rules can be formulated conditionally given the Dirichlet limit variables as independent sampling from the limit proportions (cf. Lemma 5). Furthermore, we can deduce edge selection rules for $(\mathcal{R}_k, \mu_k)$ that are analogous to (i)$_{\mathrm{rec}}$ and (ii)$_{\mathrm{rec}}$ and indeed (i) and (ii), for general $(\alpha, \theta)$.

COROLLARY 16. *Let $\theta > 0$. Then $((\mathcal{R}_k, \mu_k), k \geq 1)$ is an inhomogeneous Markov chain starting from an $(\alpha, \theta)$-string of beads $(\mathcal{R}_1, \mu_1) = (E_{\{1\}}, \mu_1)$, with transition rules, as follows:*

(i)$^{\mathcal{R}}$  *Given* $(\mathcal{R}_k, \mu_k)$, *assign weight* $\mu_k(E_B)$ *to the edge in* $(\mathcal{R}_k, \mu_k)$ *labeled* $B$, $B \in T_k$.

(ii)$^{\mathcal{R}}$  *Select* $B_k \in T_k$ *at random with probabilities proportional to the weights. Select a bead* $(J_k, m_k)$, *where* $J_k = (B_k, s_k) \in E_{B_k}$ *and* $m_k = \mu_k(\{J_k\})$ *as in Proposition 10 using* $(\alpha, \theta)$-*selection if* $\#B_k = 1$ *and* $(\alpha, \alpha)$-*selection if* $\#B_k \geq 2$ *on the string of beads* $(\widetilde{E}_B, \widetilde{\mu}_k^{E_B})$ *associated to* $(E_B, \mu_k|_{E_B})$ *by shifting and scaling.*

*To create* $\mathcal{R}_{k+1}$ *from* $\mathcal{R}_k$, *remove from* $\mathcal{R}_k$ *bead* $(J_k, m_k)$ *and attach in* $J_k$ *the* $m_k$-*scaled and* $s_k$-*shifted image* $(I_{k+1}, \mu^{I_{k+1}})$ *of an independent* $(\alpha, \theta)$-*string of beads* $(\widetilde{I}_{k+1}, \widetilde{\mu}^{I_{k+1}})$. *Relabel to include* $k+1$ *so as to obtain* $\mathcal{R}_{k+1}$ *in canonical form (22):*

$$\mathcal{R}_{k+1} = \{(A \cup \{k+1\}, s) : (A, s) \in \mathcal{R}_k, s \leq s_k, B_k \subset A\}$$
$$\cup\, I_{k+1} \cup \{(A, s) \in \mathcal{R}_k : s > s_k \text{ or } B_k \not\subset A\}$$
$$\mu_{k+1}(C) = \mu_k(\{(A, s) \in \mathcal{R}_k \setminus \{J_k\} : (A \cup \{k+1\}, s) \in C\})$$
$$+\, \mu^{I_{k+1}}(C \cap I_{k+1}) + \mu_k(\{C \cap (\mathcal{R}_k \setminus \{J_k\})\}).$$

PROOF.  $((\mathcal{R}_k, \mu_k), k \geq 1)$ is an inhomogeneous Markov chain because $(\mathcal{R}_{k+1}, \mu_{k+1})$ fully determines $(\mathcal{R}_k, \mu_k), \ldots, (\mathcal{R}_1, \mu_1)$. To identify the transition rules, fix $k \geq 1$. The proof is by induction on the steps in the recursive growth rules. The induction step consists of proving the recursive version of the growth rules (i)$^{\mathcal{R}}$ and (ii)$^{\mathcal{R}}$:

(i)$^{\mathcal{R}}_{\mathrm{rec}}$  Given $(\mathcal{R}_k, \mu_k)$ with first split $\{B, [k] \setminus B\}$, with $1 \in [k] \setminus B$ and $\#B = m$, assign weights $(\mu_k(E_{[k]}), \mu_k(\mathcal{S}_k^B), \mu_k(\mathcal{S}_k^{[k] \setminus B}))$ to the three components, that is, the trunk and the two subtrees above the first branch point.

(ii)$^{\mathcal{R}}_{\mathrm{rec}}$  Select a component at random with probabilities proportional to the weights. If a subtree with two or more leaves was selected, recursively apply the weighting procedure (i)$^{\mathcal{R}}_{\mathrm{rec}}$ to the selected subtree. Otherwise, denoting the selected edge or the unique edge in the selected subtree by $E_{B_k}$, select a bead $(J_k, m_k)$, where $J_k = (B_k, s_k) \in E_{B_k}$ and $m_k = \mu_k(\{J_k\})$ as in Proposition 10 using $(\alpha, \theta)$-selection if $\#B_k = 1$ and $(\alpha, \alpha)$-selection if $\#B_k \geq 2$ on the string of beads $(\widetilde{E}_{B_k}, \widetilde{\mu}_k^{E_{B_k}})$ associated with $(E_{B_k}, \mu_k|_{E_{B_k}})$ by shifting and scaling.

To prove that this recursive scheme produces the same distributions as the limiting procedure in Proposition 14(a) that defines $(\mathcal{R}_k, \mathcal{R}_{k+1})$, we study the independence properties in the proof of Proposition 14. The urn scheme

$$(\alpha + H_1^{(n)}, m - \alpha + H_2^{(n)}, k - m - 1 + \theta + H_3^{(n)}), \qquad n \geq k$$

starting from $H^{(k)} = (H_1^{(k)}, H_2^{(k)}, H_3^{(k)}) = (0, 0, 0)$ interacts with the growth of edges and mass measures on the subtrees only by setting the number of steps, so that by stage $n$, this growth will have exhibited $H_1^{(n)}$ steps according to the rules of ordered CRP and $H_2^{(n)}$ and $H_3^{(n)}$ steps, respectively, according to the recursive growth rules for the subtrees, irrespective of $(H^{(i)}, k \leq i < n)$. As $n \to \infty$, we obtain independence of three components $C_1, C_2, C_3$, the $(\alpha, \alpha)$-string of beads $C_1 = (\widetilde{E}_{[k]}, \widetilde{\mu}_k^{E_{[k]}})$ and the relabeled, scaled and shifted subtrees $C_2 = (\widetilde{\mathcal{S}}_k^B, \widetilde{\mu}_k^B)$ and $C_3 = (\widetilde{\mathcal{S}}_k^{[k]\backslash B}, \widetilde{\mu}_k^{[k]\backslash B})$ from the sigma-algebra $\mathcal{H}$ generated by $((H_1^{(n)}, H_2^{(n)}, H_3^{(n)}), n \geq k)$.

On the other hand, if $H_j^{(k+1)} = 1$, then leaf $k + 1$ is inserted in the $j$th component, $j = 1, 2, 3$, so this selection is $\mathcal{H}$-measurable and hence independent of $(C_1, C_2, C_3)$. Standard results on urn schemes (Lemma 5) yield that

$$\mathbb{P}(H_j^{(k+1)} = 1 | (\mathcal{R}_k, \mu_k)) = \mathbb{P}(H_j^{(k+1)} = 1 | (H_1, H_2, H_3)) = H_j \qquad \text{a.s.}$$

Inductively, this argument shows that the conditional probability given $(\mathcal{R}_k, \mu_k)$ of inserting $k + 1$ at edge $E_B$ of $\mathcal{R}_k$ is $\mu_k(E_B)$ a.s. and that, conditionally given this edge selection, the growth on that edge follows a CRP, when restricted to insertions to that edge. In particular, the bead selection is done according to Proposition 10, with parameters $(\alpha, \theta)$ if $\#B = 1$ and $(\alpha, \alpha)$ if $\#B \geq 2$; cf. Corollary 15. The insertion rule creates $E_{\{k+1\}}$ with distribution as identified in Corollary 15.  □

If $E_{B_k}$ is an internal edge, the $\mathrm{PD}(\alpha, \alpha)$ composition structure is strongly sampling consistent and, in fact, we select a new junction point $J_k$ with weights proportional to $\mu_k$ restricted to $E_{B_k}$.

For $\theta = 0$, the discussion before Proposition 10 shows that the bead selection in an $(\alpha, 0)$-string of beads always selects the last bead *at* the leaf. Crushing this bead creates a new string of beads but does not split the string the bead was selected from hence creating a degenerate subtree, which should contain the leaf edge leading to the smallest label, say, 1, for simplicity noting that this occurs recursively for all other labels also, but this edge has zero length and, in particular, no more beads. If we use the canonical representation (22), there will be no point $(\{1\}, s)$, $s \geq 0$, in $\mathcal{R}_k$, $k \geq 2$, and the "leaf" 1 is actually equal to $J_1$, a pseudo-branch point whose removal creates only two connected components. Below this point, 1 is in the label set, above it, 1 is removed from the label set.

3.6. *Moment calculations for lengths and masses.* Focusing particularly on the case $k = 2$ and $\theta = 1 - \alpha$, denote by $J_1 = (\{1, 2\}, r_{\{1,2\}})$ the branch point and by $\Sigma_i = (\{i\}, r_{\{i\}})$, $i = 1, 2$, the leaves. Then the joint distribution of lengths

$$(25) \qquad d_2(\rho, J_1), \qquad d_2(J_1, \Sigma_1), \qquad d_2(J_1, \Sigma_2)$$

was described already in [18], Proposition 18. These are dictated by the asymptotics of urn schemes embedded in the $(\alpha, 1-\alpha)$-tree growth process. In the previous subsection, we described these branch lengths jointly with the masses

$$(26) \qquad \mu_2([[\rho, J_1]]), \qquad \mu_2([[J_1, \Sigma_1]]), \qquad \mu_2([[J_1, \Sigma_2]])$$

and the restrictions of $\mu_2$ to the three branches. In the $(\alpha, 1-\alpha)$ case, Proposition 14(b) identifies the joint distribution of the sextuple (25) and (26) in terms of the Dirichlet$(\alpha, 1-\alpha, 1-\alpha)$ distribution of masses (26), and

$$(27) \qquad \begin{aligned} & d_2(\rho, J_1) = \mu_2([[\rho, J_1]])^\alpha S_0; \qquad d_2(J_1, \Sigma_1) = \mu_2([[J_1, \Sigma_1]])^\alpha S_1; \\ & d_2(J_1, \Sigma_2) = \mu_2([[J_1, \Sigma_2]])^\alpha S_2; \end{aligned}$$

where the $S_0$, $S_1$ and $S_2$ are independent $\alpha$-diversities (or local times) associated with $(\alpha, \theta)$ interval partitions with parameters $\theta = \alpha$, $\theta = 1-\alpha$ and $\theta = 1-\alpha$, respectively. It could be checked by a joint moment computation that this is consistent with the alternative description of the lengths without the masses which was provided in [18], Proposition 18:

$$(28) \qquad \begin{aligned} & d_2(\rho, J_1) = D_0 \lambda(\mathcal{R}_2); \qquad d_2(J_1, \Sigma_1) = D_1 \lambda(\mathcal{R}_2); \\ & d_2(J_1, \Sigma_2) = D_2 \lambda(\mathcal{R}_2); \end{aligned}$$

where $\lambda(\mathcal{R}_2)$ denotes the total length of $\mathcal{R}_2$ and $(D_0, D_1, D_2)$ has a Dirichlet$(1, (1-\alpha)/\alpha, (1-\alpha)/\alpha)$ distribution, independent of $\lambda(\mathcal{R}_2)$ is distributed as the $\alpha$-diversity of an $(\alpha, 2-\alpha)$ interval partition. To illustrate, the first description (27) gives

$$\mathbb{E}(d_2(\rho, J_1)^s) = \frac{B(\alpha + \alpha s, 2 - 2\alpha)}{B(\alpha, 2 - 2\alpha)} \frac{\Gamma(\alpha+1)\Gamma(s+2)}{\Gamma(2)\Gamma(\alpha + s\alpha + 1)} = \frac{\Gamma(s+1)\Gamma(2-\alpha)}{\Gamma(2 + s\alpha - \alpha)},$$

whereas the second description (28) gives

$$\mathbb{E}(d_2(\rho, J_1)^s) = \frac{B(1+s, 2/\alpha - 2)}{B(1, 2/\alpha - 2)} \frac{\Gamma(3-\alpha)\Gamma(2/\alpha + s)}{\Gamma(2/\alpha)\Gamma(3 + s\alpha - \alpha)} = \frac{\Gamma(s+1)\Gamma(2-\alpha)}{\Gamma(2 + s\alpha - \alpha)}.$$

The above discussion, together with the location of masses along the arms according to appropriate regenerative PD$(\alpha, \theta)$ distributions, with masses located at local times, fully determines the law of $(\mathcal{R}_2, \mu_2)$. What remains to be seen is how $(\mathcal{R}_2, \mu_2)$ can be embedded in the CRT.

**4. Embedding in continuum fragmentation trees.** Throughout this section we assume $0 < \alpha < 1$, since there are no CRTs (in the sense of the next subsection) associated with $\alpha = 0$ and $\alpha = 1$ (cf. the discussion at the end of Section 3.2).

4.1. *Continuum fragmentation trees.* We defined weighted $\mathbb{R}$-trees in Section 3.3. Let us follow Evans and Winter [9] to introduce a notion of convergence on the space $\mathbb{T}^{\mathrm{wt}}$ of weight-preserving isometry classes of weighted $\mathbb{R}$-trees. Here, two weighted $\mathbb{R}$-trees $(\mathcal{R}, \nu)$ and $(\mathcal{T}, \mu)$ are called weight-preserving isometric if there exists an isometry $i : \mathcal{R} \to \mathcal{T}$ with $i_* \nu = \mu$ the push-forward of measure $\nu$ under the isometry. Informally, the notion of convergence consists of weak convergence of probability measures and Gromov–Hausdorff convergence of the underlying tree spaces. See also Evans et al. [8] for Gromov–Hausdorff convergence of unweighted $\mathbb{R}$-trees and Greven, Pfaffelhuber and Winter [16] for an alternative type of convergence for weighted $\mathbb{R}$-trees.

More specifically, it is shown in [9] that the distance function

$$\Delta_{\mathrm{GH^{wt}}}((\mathcal{R}, \nu), (\mathcal{T}, \mu))$$
$$= \inf\{\varepsilon > 0 : \exists_{f \in F^\varepsilon_{\mathcal{R}, \mathcal{T}}, g \in F^\varepsilon_{\mathcal{T}, \mathcal{R}}} \, d_P(f_* \nu, \mu) \le \varepsilon \text{ and } d_P(\nu, g_* \mu) \le \varepsilon\}$$

gives rise to a Polish topology on $\mathbb{T}^{\mathrm{wt}}$ (although $\Delta_{\mathrm{GH^{wt}}}$ is not itself a metric), where

$$F^\varepsilon_{\mathcal{R}, \mathcal{T}} = \Big\{ f : \mathcal{R} \to \mathcal{T} : \sup_{x, x' \in \mathcal{R}} |d_{\mathcal{R}}(x, x') - d_{\mathcal{T}}(f(x), f(x'))| \le \varepsilon \Big\}$$

set of $\varepsilon$-isometries,

$$d_P(\mu, \mu') = \inf\{\varepsilon : \forall_{C \subset \mathcal{T} \text{ closed}} \, \mu(C) \le \mu'(\{x \in \mathcal{T} : d(x, C) \le \varepsilon\}) + \varepsilon\}$$

Prohorov distance.

Note that convergence of the form (23) for strings of beads and, based on this, (24) for sequences of weighted discrete trees with edge lengths and constant combinatorial shape imply convergence in the sense defined here. However, this notion of convergence also allows convergence to trees with more complicated branching structure such as continuum fragmentation trees.

We will further use this notion of convergence to establish projective limits of subsets of a CRT, where the measures on the subsets are just projections of the CRT mass measure. The following elementary lemma will be useful.

LEMMA 17. *Let $\mathcal{R} \subset \mathcal{T}$ be two $\mathbb{R}$-trees, $\mu$ a measure on $\mathcal{T}$ and $\nu = \pi_* \mu$ the push-forward under the projection map $\pi : \mathcal{T} \to \mathcal{R}$. Then*

$$\Delta_{\mathrm{GH^{wt}}}((\mathcal{R}, \nu), (\mathcal{T}, \mu)) \le d_{\mathrm{Haus}(\mathcal{T})}(\mathcal{R}, \mathcal{T})$$

*for the Hausdorff distance $d_{\mathrm{Haus}(\mathcal{T})}$ on compact subsets of $\mathcal{T}$.*

PROOF. Just consider the projection map $g = \pi$ and the inclusion map $f : \mathcal{R} \to \mathcal{T}$, then for $\varepsilon = d_{\mathrm{Haus}(\mathcal{T})}(\mathcal{R}, \mathcal{T})$, we have $f \in F^{\varepsilon}_{\mathcal{R}, \mathcal{T}}, g \in F^{\varepsilon}_{\mathcal{T}, \mathcal{R}}, d_P(\nu, g_* \mu) = 0$ and $d_P(f_* \nu, \mu) = d_P(\nu, \mu) \le \varepsilon$. $\quad\square$

A random weighted rooted *binary* $\mathbb{R}$-tree $(\mathcal{T}, d, \rho, \mu)$ is called a *binary fragmentation CRT* of index $\gamma > 0$, if

- $\mu$ is nonatomic a.s. assigning positive weight to the subtrees $\mathcal{T}_\sigma = \{\sigma' \in \mathcal{T} : \sigma \in [[\rho, \sigma']]\}$ for all nonleaf $\sigma \in \mathcal{T}$, and zero weight to all branches $[[\rho, \sigma]]$, for all $\sigma \in \mathcal{T}$, and
- for all $t \ge 0$ the connected components $(\mathcal{T}_i^t, i \ge 1)$ of $\{\sigma \in \mathcal{T} : d(\rho, \sigma) > t\}$, completed by a root vertex $\rho_i$, are such that given $(\mu(\mathcal{T}_i^t), i \ge 1) = (m_i, i \ge 1)$ for some $m_1 \ge m_2 \ge \cdots \ge 0$, the trees

$$(\mathcal{T}_i^t, m_i^{-\gamma} d|_{\mathcal{T}_i^t}, \rho_i, m_i^{-1} \mu|_{\mathcal{T}_i^t}), \qquad i \ge 1,$$

are like independent identically distributed isometric copies of $(\mathcal{T}, d, \rho, \mu)$.

Haas and Miermont [17] and Bertoin [3] observed the following. Given $(\mathcal{T}, d, \mu)$, let $\Sigma^*$ be a random point in $\mathcal{T}$ chosen according to $\mu$, and define the mass of the *tagged* subtree above $t$ as

$$S_t^* = \begin{cases} \mu(\mathcal{T}_i^t), & \text{if } \Sigma^* \in \mathcal{T}_i^t \text{ for some } i \ge 1, \\ 0, & \text{otherwise.} \end{cases}$$

Then $(S_t^*, t \ge 0)$ is a decreasing self-similar Markov process in $[0, 1]$ starting from $S_0^* = 1$ and attaining $S_t^* = 0$ in finite time, which can be expressed as

$$S_t^* = \exp\{-\xi^*_{T(t)}\} \qquad \text{where } T(t) = \inf\left\{ u \ge 0 : \int_0^u \exp\{-\gamma \xi^*_r\} \, dr > t \right\}$$

and $\xi^*$ is a subordinator, called the *spinal subordinator*, with Laplace exponent

$$\Phi(s) = \int_{(0,\infty)} (1 - e^{-sx}) \Lambda^*(dx)$$

for some Lévy measure $\Lambda^*$ on $(0, \infty)$ with $\int_{(0,\infty)} (1 \wedge x) \Lambda^*(dx) < \infty$ that characterizes the distribution of the *binary* fragmentation CRT. A jump $\Delta \xi^*_{T(t)} = x$ corresponds to a change of mass $S_t^* = S_{t-}^* e^{-x}$ by a factor of $e^{-x}$ at height $t$, so consider the push-forward $\tilde{\Lambda}^*(du)$ of $\Lambda^*$ via the transformation $u = e^{-x}$. It will be assumed in the following discussion that $\Lambda^*(dx) = \lambda^*(x) \, dx$ for some density function $\lambda^*(x)$, so that $\tilde{\Lambda}^*(du) = u f^*(u) \, du$ for some density function $f^*$ on $(0, 1)$ which is related to $\lambda^*$ by

$$(29) \qquad f^*(u) = u^{-2} \lambda^*(-\log u) \qquad (0 < u < 1).$$

The introduction of the size-biasing factor $u$ is done since the normal parameterization of fragmentation trees is by their *dislocation measure*

$$\nu(du) = 1_{\{u \geq 1/2\}} f^*(u) \, du.$$

The size-biasing factor $u$ then arises because in our context of binary fragmentations, $f^*$ is necessarily *symmetric*, meaning $f^*(u) = f^*(1-u)$, and given a mass split $(u, 1-u)$ with $u < 1-u$, the mass of the randomly tagged fragment is multiplied by $u$ with probability $u$ and by $1-u$ with probability $1-u$, but then the total rate for a ranked split $(u, 1-u)$ with $u \geq 1/2$ is again $u f^*(u) + (1-u) f^*(1-u) = f^*(u)$.

Because $\xi^*$ is a subordinator, $\{1 - \exp(-\xi_t^*), t \geq 0\}^{\mathrm{cl}}$ is a regenerative interval partition in the sense of Section 2.1.

PROPOSITION 18 (Spinal decomposition [4, 19]). *Consider a fragmentation CRT $(\mathcal{T}, d, \mu)$ and a random leaf $\Sigma^* \in \mathcal{T}$ whose distribution given $(\mathcal{T}, d, \mu)$ is $\mu$. Then the spinal decomposition theorem holds for the spine $[[\rho, \Sigma^*]]$ in the following sense. Consider the connected components $(\mathcal{T}_i, i \in I)$ of $\mathcal{T} \setminus [[\rho, \Sigma^*]]$, each completed by a root vertex $\rho_i$. Denote by $\mu^*$ the random discrete distribution on $[[\rho, \Sigma^*]]$ obtained by assigning mass $m_i = \mu(\mathcal{T}_i)$ to the branch point base point of $\mathcal{T}_i$ on $[[\rho, \Sigma^*]]$. Then given the string of beads $([[\rho, \Sigma^*]], \mu^*)$, the trees*

$$(\mathcal{T}_i, m_i^{-\gamma} d|_{\mathcal{T}_i}, \rho_i, m_i^{-1} \mu|_{\mathcal{T}_i}), \qquad i \in I,$$

*are independent identically distributed isometric copies of $(\mathcal{T}, d, \mu)$.*

### 4.2. $(\alpha, \theta)$-dislocation measures and switching probabilities.

From Proposition 14 we have $(\alpha, \theta)$-trees $(\mathcal{R}_k, \mu_k)$ which are based on weakly sampling consistent regenerative Poisson–Dirichlet compositions. We can compare this with sampling $k$ leaves $\Sigma_1^*, \ldots, \Sigma_k^*$ according to $\mu$ in a CRT $(\mathcal{T}, \mu)$ giving rise to reduced fragmentation trees

$$\mathcal{R}_k^* = \bigcup_{j=1}^{k} [[\rho, \Sigma_j^*]], \qquad \mu_k^* = \pi_*^{\mathcal{R}_k^*} \mu,$$

which can be thought of as being based on strongly sampling consistent regenerative compositions that are *not* of Poisson–Dirichlet type [by Proposition 6(ii), the unique *regenerative* Poisson–Dirichlet interval partition is not strongly sampling consistent unless $\alpha = \theta = 1/2$].

Let us first compute the appropriate dislocation measure for $(\mathcal{T}, \mu)$. We have a (spinal-to-be) subordinator $\xi$ with Laplace exponent $c\Phi_{\alpha, \theta}$ given by (10), with Lévy measure

$$\Lambda_{\alpha, \theta}(dx) = \lambda_{\alpha, \theta}(x) \, dx,$$

where $\lambda_{\alpha,\theta}(x) = c\alpha(1 - e^{-x})^{-\alpha-1}e^{-(\theta+1)x} + c\theta e^{-\theta x}(1 - e^{-x})^{-\alpha}$. Here, $c$ is a constant that was irrelevant in the context of Section 2, but that we will choose appropriately here. In analogy with (29), we can compute the intensity $uf_{\alpha,\theta}(u)$ of $(e^{-\Delta\xi_t})$ as

$$f_{\alpha,\theta}(u) = u^{-2}\lambda_{\alpha,\theta}(-\log(u)) = c\alpha(1-u)^{-\alpha-1}u^{\theta-1} + c\theta u^{\theta-2}(1-u)^{-\alpha},$$

which is nonsymmetric and, for a split $(u, 1-u)$ with $u \geq 1/2$, gives a rate of

$$\begin{aligned}
f^\circ_{\alpha,\theta}(u) &= uf_{\alpha,\theta}(u) + (1-u)f_{\alpha,\theta}(1-u) \\
&= c\alpha(u^\theta(1-u)^{-\alpha-1} + (1-u)^\theta u^{-\alpha-1}) \\
&\quad + c\theta(u^{\theta-1}(1-u)^{-\alpha} + (1-u)^{\theta-1}u^{-\alpha}).
\end{aligned}$$

We now check that the choice $c = 1/\Gamma(1-\alpha)$ is such that

$$\nu_{\alpha,\theta}([1/2, 1-\varepsilon]) := \int_{1/2}^{1-\varepsilon} f^\circ_{\alpha,\theta}(u)\,du \sim \frac{\varepsilon^{-\alpha}}{\Gamma(1-\alpha)},$$

which is the condition established in [18] to obtain the associated CRT as limit in discrete approximations scaled by $n^\alpha$ as in Proposition 14, but in the weakly sampling consistent case.

We can now compare subordinators $\xi$ with Lévy measure $\Lambda_{\alpha,\theta}$ and the spinal subordinator $\xi^*$ in a CRT $(\mathcal{T}, \mu)$ with dislocation density $f^\circ_{\alpha,\theta}$. Recall also [14, 18] that the regenerative interval partition associated with the spinal subordinator $\xi^*$ admits a natural local time process as in (13) and (14), which is such that the spinal string of beads $([[\rho, \Sigma^*]], \mu^*)$ is of the form

$$\begin{aligned}
&d(\rho, \Sigma^*) = L^*(1) \quad \text{and} \\
&\mu^*(\{g_{\rho,\Sigma^*}(L^*(1 - e^{-\xi^*_\tau}))\}) = e^{-\xi^*_{\tau-*}} - e^{-\xi^*_\tau} = -\Delta e^{-\xi^*_\tau}.
\end{aligned} \tag{30}$$

In particular, we can identify the height $L^*(1 - e^{-\xi^*_\tau})$ in the tree of an atom of $\mu^*$ that corresponds to a jump of $\xi^*$ at time $\tau$.

For the proof of Theorem 4, we will embed $(\mathcal{R}_k, \mu_k)$, $k \geq 1$, in the CRT $(\mathcal{T}, \mu)$ of index $\alpha$ and with dislocation density $f^\circ_{\alpha,\theta}$. This involves solving several problems:

- How do we embed $(\mathcal{R}_1, \mu_1)$ in $(\mathcal{T}, \mu)$? Can we make leaf $\Sigma_1$ of $\mathcal{R}_1$ close to leaf $\Sigma_1^*$ of $\mathcal{R}_1^*$ by having their spines coincide initially? Part of the problem is then to identify the point where the spines separate.
- Can we iterate the procedure by following the exchangeable leaf with the smallest label $\Sigma_{n_{i,1}}^*$ off the spine of $\mathcal{R}_1^*$, and pass to a limit $i \to \infty$ to identify $\mathcal{R}_1$ as a subset of $\mathcal{T}$?
- Once we have $(\mathcal{R}_1, \mu_1)$, how do we find the point where the spine of leaf $\Sigma_2$ leaves the spine of leaf $\Sigma_1$?

• Can we iterate this to embed all $(\mathcal{R}_k, \mu_k)$ in $(\mathcal{T}, \mu)$?

Outside a CRT, we solved the third bullet point in Proposition 10 and obtained a coin-tossing representation in the sense that we climb up the spine tossing a coin for each of the (infinite number of subtree) masses and stopping the first time we see heads. The heads probability depends on the relative remaining mass $u$ after a split and can be given as a *switching probability* (away from relative size $u$ to relative size $1 - u$)

$$p(u) = \frac{(1 - u)\theta}{(1 - u)\theta + u\alpha} \qquad \text{where } u = \exp(-\Delta\xi_t).$$

See also Corollary 16 for the iteration for $k \geq 2$. Although we endeavor to embed $(\mathcal{R}_1, \mu_1)$ in $(\mathcal{T}, \mu)$, it is instructive to first try to embed $(\mathcal{R}_1^*, \mu_1^*)$ in $(\mathcal{R}_n, \mu_n)$. Assuming for a moment that $\mathcal{R}_n \subset \mathcal{T}$ and $\mu_n = \pi_*^{\mathcal{R}_n}\mu$, then $\Sigma_1^*$ as a pick from $\mu$ is projected onto $\Sigma_{1,n}^* = \pi^{\mathcal{R}_n}(\Sigma_1^*)$, a pick from $\mu_n$.

LEMMA 19. (a) *Given* $(\mathcal{R}_1, \mu_1)$, *a pick* $\Sigma_{1,1}^*$ *from* $\mu_1$ *is obtained by switching probabilities* $p^*(u) = 1 - u$: *given* $(\mathcal{R}_1, \mu_1)$ *is associated with a spinal subordinator* $\xi$, *the conditional probability that* $\Sigma_{1,1}^*$ *falls into the block* $(1 - e^{-\xi_{t-}}, 1 - e^{-\xi_t})$ *of the associated interval partition is*

$$p^*(e^{-\Delta\xi_t}) \prod_{s < t}(1 - p^*(e^{-\Delta\xi_s})).$$

(b) *Let* $\theta > 0$. *Denote the switching time in* (a) *by* $\tau$. *Given* $(\mathcal{R}_1^*, \mu_1^*)$ *and a measurable switching probability function* $(\widehat{p}(u), 0 \leq u \leq 1)$ *with associated switching time* $\widehat{\tau}$, *we obtain*

$$(31) \qquad (\xi_t, 0 \leq t < \tau) \stackrel{d}{=} (\xi_t^*, 0 \leq t < \widehat{\tau})$$

*if and only if*

$$(32) \qquad \widehat{p}(u) = \frac{(1 - u)f_{\alpha,\theta}(1 - u)}{f_{\alpha,\theta}^*(u)} \qquad \text{for almost all } 0 \leq u \leq 1.$$

PROOF. For (a) just note that

$$(1 - e^{-\xi_t}) - (1 - e^{-\xi_{t-}}) = e^{-\xi_{t-}}(1 - e^{-\Delta\xi_t}) = (1 - e^{-\Delta\xi_t}) \prod_{s < t} e^{-\Delta\xi_s}.$$

For (b), note that the killed subordinator $(\xi_t, 0 \leq t < \tau)$ can be described in terms of two independent Poisson point processes of points $e^{-\Delta\xi_t}$ with tails coin toss at intensity measure $(1 - p^*(u))u f_{\alpha,\theta}(u)\,du$ and jumps with heads coin toss at total intensity $\int_{(0,1)} p^*(u)u f_{\alpha,\theta}(u)\,du$.

Similarly, the killed subordinator $(\xi_t^*, 0 \leq t < \widehat{\tau})$ has the same description with tails intensity measure $(1 - \widehat{p}(u))u f_{\alpha,\theta}^*(u)\,du$ and heads intensity $\int_{(0,1)} \widehat{p}(u) u f^*_{\alpha,\theta}(u) \, du$. It is an elementary computation to show that the tails intensity measures are equal if and only if $\widehat{p}(u)$ satisfies (32) and that then also the heads intensities coincide. $\quad\square$

For $\theta = 0$, the subordinator $\xi$ with Laplace exponent (10) has an infinite jump $\Delta \xi_{\mathbf{e}} = \infty$ at an exponential time $\mathbf{e}$ with parameter $1/\Gamma(1-\alpha)$, while $\xi^*$ does not. The calculation in the proof is still true, except that the possibility of an infinite jump was ignored. Consequently, for (31) to hold, $\widehat{\tau}$ must be replaced by $\widehat{\tau} \wedge \mathbf{e}$ for an independent exponential time $\mathbf{e}$ with parameter $1/\Gamma(1-\alpha)$, that is,

$$(33) \qquad (\xi_t, 0 \le t < \tau) \overset{d}{=} (\xi^*_t, 0 \le t < \widehat{\tau} \wedge \mathbf{e}).$$

Note that $\mathbf{e}$ is not a jump time of $\xi^*$.

4.3. *Embedding $(\mathcal{R}_k, \mu_k)$ and the proof of Theorem 4.* We now carry out the program outlined in the previous subsection and iterate the embedding started in Lemma 19 to construct an unkilled Poisson point process $(F_t, t \ge 0)$ and then $(\mathcal{R}_1, \mu_1)$:

- Let $(\mathcal{T}, d, \rho, \mu)$ be an $\alpha$-self-similar fragmentation CRT with dislocation density $f^\circ_{\alpha,\theta}$.
- Define $(\mathcal{T}^{(1)}, d^{(1)}, \rho^{(1)}, \mu^{(1)}) := (\mathcal{T}, d, \rho, \mu)$ and consider the spinal subordinator $\xi^{*(1)}$ of a random point $\Sigma^{*(1)}_1$ sampled from $\mu^{(1)}$ in $\mathcal{T}^{(1)}$. Perform the construction of Lemma 19(b) and denote by $\tau^{(1)}$ the associated switching time, also put $\tau^{(0)} = 0$. Define

  $$F_t = \exp(-\Delta \xi^{*(1)}_t) \qquad \text{for } 0 \le t < \tau^{(1)}, \qquad F_{\tau^{(1)}} = 1 - \exp(-\Delta \xi^{*(1)}_{\tau^{(1)}}).$$

  For $\theta = 0$, when $\tau^{(1)} = \widehat{\tau} \wedge \mathbf{e}$ in (33), terminate the construction if $\tau^{(1)} = \mathbf{e}$.
- For $i \ge 1$, denote by $(L^{*(i)}(u), 0 \le u \le 1)$ the local time process associated with the interval partition $\{1 - e^{-\xi^{*(i)}_t}, t \ge 0\}^{\mathrm{cl}}$ and by

  $$\rho^{(i+1)} = g_{\rho^{(i)}, \Sigma^{*(i)}}(1 - L^{*(i)}(\exp(-\xi^{*(i)}_{\tau^{(i)}})))$$

  the junction point; cf. (30). Define

  $$\mathcal{T}^{(i+1)} = \{\sigma \in \mathcal{T}^{(i)} : [[\rho^{(i)}, \sigma]] \cap [[\rho^{(i)}, \Sigma^{*(i)}_1]] = [[\rho^{(i)}, \rho^{(i+1)}]]\},$$

  $$d^{(i+1)} = (1 - \exp(-\xi^{*(i)}_{\tau^{(i)} - \tau^{(i-1)}}))^{-\alpha} d^{(i)}|_{\mathcal{T}^{(i+1)}},$$

  $$\mu^{(i+1)} = (1 - \exp(-\xi^{*(i)}_{\tau^{(i)} - \tau^{(i-1)}}))^{-1} \mu^{(i)}|_{\mathcal{T}^{(i+1)}}.$$

Then consider the spinal subordinator $\xi^{*(i+1)}$ of $\Sigma_1^{*(i+1)} \sim \mu^{(i+1)}$ in $\mathcal{T}^{(i+1)}$. Perform the construction of Lemma 19(b) and denote by $\tau^{(i+1)}$ the associated switching time. Define

$$F_{\tau^{(i)}+t} = \exp(-\Delta\xi_t^{*(i+1)}) \qquad \text{for } 0 \le t < \tau^{(i+1)} - \tau^{(i)},$$

$$F_{\tau^{(i+1)}} = 1 - \exp(-\Delta\xi_{\tau^{(i+1)}-\tau^{(i)}}^{*(i+1)}).$$

PROPOSITION 20. (a) *For* $\theta > 0$, *the process* $(F_t, t \ge 0)$ *is a Poisson point process with intensity measure* $uf_{\alpha,\theta}(u)$ *(and cemetery state* 1*). The subspace* $[[\rho, \Sigma_1[[:= \bigcup_{i\ge1}[[\rho, \rho^{(i)}]]$ *is such that* $\Sigma_1 \in \mathcal{T}$ *is a leaf a.s., and* $([[\rho, \Sigma_1]], \pi_*^{[[\rho,\Sigma_1]]}\mu)$ *is a weight-preserving isometric copy of* $(\mathcal{R}_1, \mu_1)$. *Furthermore, the spinal decomposition theorem holds for the spine* $[[\rho, \Sigma_1]]$ *and the connected components* $(\mathcal{T}_i, i \in I)$ *of* $\mathcal{T} \setminus [[\rho, \Sigma_1]]$; *cf. Proposition* 18.

(b) *For* $\theta = 0$, *the process* $(F_t, 0 \le t < \mathbf{e})$ *is a Poisson point process with intensity measure* $uf_{\alpha,\theta}(u)$ *killed at an independent exponential time* $\mathbf{e}$ *with parameter* $1/\Gamma(1-\alpha)$. *Denote* $I$ *such that* $\tau^{(I)} = \mathbf{e}$. *Then the subspace* $[[\rho, \Sigma_1]] := [[\rho, \rho^{(I)}]]$ *is such that* $\Sigma_1$ *is not a leaf a.s. The weighted space* $([[\rho, \Sigma_1]], \pi_*^{[[\rho,\Sigma_1]]}\mu)$ *is an isometric copy of* $(\mathcal{R}_1, \mu_1)$. *The spinal decomposition theorem holds for the spine* $[[\rho, \Sigma_1]]$.

PROOF. By Lemma 19, $(F_t, 0 \le t < \tau^{(1)})$ is a Poisson point process with intensity parameter $u^2 f_{\alpha,\theta}(u)\,du$ killed at an independent exponential time with parameter $\kappa = \int_{(0,1)}(1-u)uf_{\alpha,\theta}(u)\,du$, and $F_{\tau^{(1)}}$ has distribution $\kappa^{-1}\hat{p}(1-u)(1-u)f_{\alpha,\theta}^*(1-u) = \kappa^{-1}(1-u)uf_{\alpha,\theta}(u)\,du$.

For $i \ge 1$, denote by $\mathcal{G}_i = \sigma((\xi^{*(1)}, \tau^{(1)}), \ldots, (\xi^{*(i)}, \tau^{(i)}))$ the $\sigma$-algebra generated by the first $i$ spinal subordinators and their switching times. It follows easily from the definition that $\mathcal{T}^{(i+1)} \setminus \{\rho^{(i+1)}\}$ is a connected component of $\mathcal{T}^{(i)} \setminus [[\rho^{(i)}, \Sigma_1^{*(i)}]]$. By Proposition 18, the tree $(\mathcal{T}^{(i+1)}, d^{(i+1)}, \rho^{(i+1)}, \mu^{(i+1)})$ is a copy of $(\mathcal{T}^{(i)}, d^{(i)}, \rho^{(i)}, \mu^{(i)})$ that is independent of $\mathcal{G}_i$.

By induction and standard superposition results for Poisson point processes, the process $(F_t, t \ge 0)$ is a Poisson point process with intensity measure

$$u^2 f_{\alpha,\theta}(u)\,du + (1-u)uf_{\alpha,\theta}(u)\,du = uf_{\alpha,\theta}(u)\,du,$$

as claimed. In particular, the associated mass process $e^{-\xi_t} = \prod_{s\le t}F_s$ has the same distribution as the process associated with $(\mathcal{R}_1, \mu_1)$.

For $\theta > 0$, completeness of $\mathcal{T}$ implies $\Sigma_1 \in \mathcal{T}$ and $e^{-\xi_\infty} = 0$ yields that $\Sigma_1$ is a leaf, since $\mu$ would otherwise assign positive mass to the subtree $\mathcal{T}_{\Sigma_1}$ above $\Sigma_1$. For $\theta = 0$, note that $\Sigma_1 \in \mathcal{T}_{\rho^{(I)}} \setminus \{\rho^{(I)}\}$.

The spinal decomposition theorem follows by a simple induction, from a version of Proposition 18 where $\Sigma^*$ is replaced by $\rho^{(2)}$. That result is proved like Proposition 18 using partition-valued fragmentation processes and stopping lines; see [4, 19]. □

The interval $[[\rho, \Sigma_1]]$ has length

$$(34) \qquad d(\rho, \Sigma_1) = \Gamma(1-\alpha) \int_0^\infty \exp(-\alpha \xi_t) \, dt = L(1),$$

whereas the interval $[[\rho, \Sigma_1^*]]$ has length

$$(35) \qquad d(\rho, \Sigma_1^*) = \Gamma(1-\alpha) \int_0^\infty \exp(-\alpha \xi_t^{*(1)}) \, dt = L^{*(1)}(1).$$

We have joined these two intervals at a junction point $J_{1,1*} = \rho^{(2)}$ at distance

$$(36) \qquad \begin{aligned} d(\rho, J_{1,1*}) &= \Gamma(1-\alpha) \int_0^{\tau^{(1)}} \exp(-\alpha \xi_t) \, dt \\ &= \Gamma(1-\alpha) \int_0^{\tau^{(1)}} \exp(-\alpha \xi_t^{*(1)}) \, dt, \end{aligned}$$

where $\tau^{(1)}$ is the switching time for the two coupled subordinators. Now the points $\Sigma_1, \Sigma_1^*, J_{1,1*}$ have been embedded in the CRT $(\mathcal{T}, \mu)$.

So $\mathcal{R}_1$ and $\mathcal{R}_1^*$ are both embedded as paths in $(\mathcal{T}, \mu)$. Moreover, if we consider the strings of beads $(\mathcal{R}_1, \mu_1)$ and $(\mathcal{R}_1, \mu_1^*)$ associated via (30), the measures $\mu_1$ and $\mu_1^*$ are the projections onto $\mathcal{R}_1$ and $\mathcal{R}_1^*$ of the mass measure $\mu$ in the CRT $(\mathcal{T}, \mu)$. We can now check that, for $\theta = 1 - \alpha$, the random length $d(\rho, \Sigma_1)$ in (34) has the same distribution as the length $S_1$ described in [18], Proposition 18. From previous discussions, the ranked masses of $\mu_1$ have PD$(\alpha, \theta)$ distribution. The interval partition of $[0, 1]$ obtained by putting these masses in the order they appear along $\mathcal{R}_1 = [[\rho, \Sigma_1]]$ is that associated with an $(\alpha, \theta)$ regenerative composition of $[0, 1]$.

Turning to $k = 2$, we identified switching probabilities in Proposition 10 that identify the branch point for $\mathcal{R}_2$ in $\mathcal{R}_1$. As $\mathcal{R}_1$ has been embedded in $\mathcal{T}$, we identify the branch point in $\mathcal{T}$. Since the spinal decomposition theorem holds for the spine $[[\rho, \Sigma_1]]$, to embed $\Sigma_2$, we repeat in the subtree thus identified the procedure we used to embed $\Sigma_1$ in $\mathcal{T}$. In particular, this procedure also constructs the mass measure $\mu_2$ as the projection onto $\mathcal{R}_2$ of the mass measure $\mu$ on the CRT.

An inductive step from $(R_k, \mu_k)$ to $(R_{k+1}, \mu_{k+1})$ now completes the embedding and hence the proof of Theorem 4. The inductive assumption will be that $(R_k, \mu_k)$ has been embedded in the CRT with $\mu_k$ the projection of the mass measure $\mu$ of $\mathcal{T}$, along with a description of $\mu_k$ as in Proposition 14.

This establishes the following corollary to Proposition 20.

COROLLARY 21.  *Given $(\mathcal{R}_k, \mu_k)$ embedded in $(\mathcal{T}, \mu)$, proceed as in Corollary 16: first pick an edge according to the allocation of mass to edges by $\mu_k$. If the edge is an inner edge, pick $J_k$ from $\mu_k$ conditioned on that edge. If the edge is a leaf edge, pick $J_k$ instead from the atoms of $\mu_k$ on this edge according to the scheme used to pick $J_1$ from $\mathcal{R}_1$, using the obvious bijection. In either case, distribute the mass $\mu_k(\{J_k\})$ onto a new edge $[[J_k, \Sigma_{k+1}]]$ according to a scaled copy of the construction of $\mathcal{R}_1$ in Proposition 20. Then the tree $(\mathcal{R}_k \cup [[J_k, \Sigma_{k+1}]])$ with measure as described is a copy of $(\mathcal{R}_{k+1}, \mu_{k+1})$.*

PROOF OF THEOREM 4.  The embedding of $(\mathcal{R}_1, \mu_1)$ into $(\mathcal{T}, \mu)$ was given in Proposition 20. An induction based on Corollary 21 completes the embedding of $(\mathcal{R}_k, \mu_k)$, $k \geq 1$.  □

4.4. *Convergence of Markov branching trees and the proof of Theorem 3.* An attractive feature of the above construction is that by a fairly obvious extension we can construct an $\mathcal{R}_k$ spanned by a root and $\Sigma_1, \ldots, \Sigma_k$ governed by the $(\alpha, \theta)$-rules, and a leaf exchangeable $\mathcal{R}_k^*$ spanned by a root and $\Sigma_1^*, \ldots, \Sigma_k^*$, all embedded in the same CRT $(\mathcal{T}, \mu)$. Specifically, $\Sigma_{k+1}^*$ and $\Sigma_{k+1}$ will by construction project onto the same edge of $\mathcal{R}_k$.

PROPOSITION 22.  *In the above construction, $d(\Sigma_k, \Sigma_k^*) \to 0$ almost surely.*

PROOF.  We work conditionally given $(\mathcal{T}, \mu)$. Let $\theta > 0$. Let us show that, for all $\varepsilon > 0$, there a.s. exists $k_1 \geq 1$ such that all edges of $\mathcal{R}_{k_1}$ have length less than $\varepsilon/3$ and all connected components of $\mathcal{T} \setminus \mathcal{R}_{k_1}$ have diameter less than $\varepsilon/3$.

First, to fix a subtree of diameter $\varepsilon/3$, consider the connected components of

$$\{\sigma \in \mathcal{T} : \{\sigma' \in \mathcal{T}_\sigma : d(\sigma, \sigma') \geq \varepsilon/3\} = \varnothing\},$$

each completed by their root on the branches of $\mathcal{T}$. Since $\mathcal{T}$ is compact, at most finitely many components $\mathcal{T}_1, \ldots, \mathcal{T}_N$ actually attain height $\varepsilon/3$. Fix subtree $\mathcal{T}_j$ with root $R_j$, and denote its mass by $m_j$. Note that the interval partitions $\mathcal{Z}_\sigma$, $\sigma \in \mathcal{T}_j \setminus \{R_j\}$, induced by $([[\rho, \sigma]], \pi_*^{[\rho, \sigma]} \mu)$ coincide on $[0, 1 - m_j]$, and denote the components of the restricted interval partition $[0, 1 - m_j] \setminus \mathcal{Z}_\sigma = \bigcup_{i \in \mathcal{I}_j : g_i \leq 1 - m_j} (g_i, d_i)$. Now, in the notation of the proof of Proposition 10,

$$q_k := \mathbb{P}(T_j \cap (\mathcal{R}_k \setminus \mathcal{R}_{k-1}) \neq \varnothing | \mathcal{R}_{k-1}) \geq m_j \prod_{i \in \mathcal{I}_j : g_i \leq 1 - m_j} (1 - p_{g_i})  \qquad \text{a.s.}$$

is bounded below uniformly in $k$. Therefore, the step when $\mathcal{T}_j \cap \mathcal{R}_k \neq \varnothing$ is bounded by a geometric random variable, and no subtrees of height $\varepsilon/3$ can persist outside $\mathcal{R}_k$ forever, so there a.s. exists $k_0 \geq 1$ such that $\mathcal{T} \setminus \mathcal{R}_{k_0}$ has no connected components of diameter exceeding $\varepsilon/3$.

Second, fix an edge of $\mathcal{R}_{k_0}$ of length exceeding $\varepsilon/3$. There are at most $2k_0 - 1$ such edges. The projected mass is an $(\alpha, \alpha)$ or $(\alpha, \theta)$-string of beads, dense on the edge. The dynamics of the growth process in Corollary 21 are such that cut points on inner edges are selected according to the mass distribution. On leaf edges, an argument as for subtrees applies. Note also that all edges added in the growth procedure after step $k_0$ are part of a subtree of diameter less than $\varepsilon/3$ and hence shorter than $\varepsilon/3$. Therefore, there a.s. exists $k_1 \geq k_0$ such that all edges in $\mathcal{R}_k$ are shorter than $\varepsilon/3$ for all $k \geq k_1$.

In particular, for all $k \geq k_1$, we deduce $d(\Sigma_{k+1}, \Sigma_{k+1}^*) < \varepsilon$ a.s., as required.

For $\theta = 0$, the arguments still apply, but some details are different. Specifically, the first time a leaf edge is picked, the atom at its top is selected and spread over a new edge, the original edge then being an internal edge and the above argument applies. Similarly, the lower bound given for $q_n$ will vanish if $R_j$ is an interior point of a leaf edge of $\mathcal{R}_{k-1}$; but we can then proceed in two steps. Specifically, we first pick this leaf edge after a geometric time, when the mass at its leaf is spread over a new edge, the original edge then being an internal edge and $\mathcal{T}_j$ is then attained after a further geometric time with parameter $m_j$. $\square$

PROOF OF THEOREM 3. The argument given in the proof of Proposition 22 also shows that $\mathcal{R}_k$ converges to $\mathcal{T}$ a.s. in the Hausdorff sense, which implies convergence of their isometry classes in the Gromov–Hausdorff sense. This proves the statement of Theorem 3 for the trees $\mathcal{R}_k$ constructed in Theorem 4, which *assumes* the existence of a CRT $(\mathcal{T}, \mu)$ on the given probability space and sufficient extra randomness to sample repeatedly from $\mu$ as needed for the construction of $\mathcal{R}_k$.

If $\mathcal{R}_k$, $k \geq 1$, are constructed from an $(\alpha, \theta)$-tree growth process as in Proposition 2, then we use the fact that the whole sequence $(\mathcal{R}_k, k \geq 1)$ has the same distribution as if it was constructed as above. Almost sure convergence in the Gromov–Hausdorff sense is a property of the distribution on $\mathbb{T}^{\mathbb{N}}$, where $\mathbb{T}$ denotes the space of isometry classes of compact real trees. We can define the limiting $\mathbb{R}$-tree $\mathcal{T}$ as the metric completion of $\bigcup_{k \geq 1} \mathcal{R}_k$, using the completeness of $\mathbb{T}$. $\square$

Another consequence is that the uniform measure on leaves of $\mathcal{R}_k$ is closely coupled to the uniform measure on leaves of $\mathcal{R}_k^*$, and hence to the mass measure $\mu$ in the CRT.

COROLLARY 23. *In the setting of Proposition 2, there exists a CRT $(\mathcal{T}, \mu)$ on the same probability space, such that following convergences hold:*

$$(\mathcal{R}_k, \mu_k) \to (\mathcal{T}, \mu) \qquad \text{in the weighted Gromov–Hausdorff sense,}$$

where $\mu_k$ is the measure identified in Proposition 14(a), and

$$(\mathcal{R}_k, \nu_k) \to (\mathcal{T}, \mu) \qquad \text{in the weighted Gromov–Hausdorff sense,}$$

where $\nu_k$ is the empirical measure on the $k$ leaves of $\mathcal{R}_k$.

PROOF. We prove this for the embedded versions of Theorem 4. Since $\mu_k$ is the projection of $\mu$ onto $\mathcal{R}_k \subset \mathcal{T}$, the first convergence is a direct consequence of Lemma 17 and the proof of Theorem 3.

For the second convergence fix $\varepsilon > 0$. Let $k_1 \geq 1$ such that $\mathcal{T} \setminus \mathcal{R}_k$ has no subtrees of diameter exceeding $\varepsilon/9$ and, hence, $d(\Sigma_k, \Sigma_k^*) < \varepsilon/3$ for $k \geq k_1$. Let $k_2 \geq 3k_1/\varepsilon$ and $k_3 \geq k_2$ such that $d_{\mathrm{Haus}(\mathcal{T})}(\mathcal{R}_k, \mathcal{T}) < \varepsilon$ for $k \geq k_3$. Then the triangular inequality for the Prohorov distance shows for $g = \pi^{\mathcal{R}_k}$ and $f : \mathcal{R}_k \to \mathcal{T}$ the inclusion map that

$$d_P(f_* \nu_k, \mu) = d_P(\nu_k, \mu) \leq d_P(\nu_k, \mu_k) + d_P(\mu_k, \mu) \leq \varepsilon$$

and

$$d_P(\nu_k, g_* \mu) \leq d_P(\nu_k, \mu_k) + d_P(\mu_k, \mu) + d_P(\mu, g_* \mu) < \varepsilon$$

for all $k \geq k_2$. This completes the proof. □

DEPARTMENT OF STATISTICS
UNIVERSITY OF CALIFORNIA AT BERKELEY
BERKELEY, CALIFORNIA 94720
USA
E-MAIL: pitman@stat.berkeley.edu

DEPARTMENT OF STATISTICS
1 SOUTH PARKS ROAD
OXFORD OX1 3TG
UNITED KINGDOM
E-MAIL: winkel@stats.ox.ac.uk